


\magnification 1200
\hsize 13.2cm
\vsize 20cm
\parskip 3pt plus 1pt
\parindent 5mm

\def\\{\hfil\break}


\font\seventeenbf=cmbx10 at 17.28pt

\font\twelvebf=cmbx10 at 12pt
\font\eightbf=cmbx8
\font\sixbf=cmbx6

\font\eighti=cmmi8
\font\sixi=cmmi6

\font\eightrm=cmr8
\font\sixrm=cmr6

\font\eightsy=cmsy8
\font\sixsy=cmsy6

\font\eightit=cmti8
\font\eighttt=cmtt8
\font\eightsl=cmsl8

\font\seventeenbsy=cmbsy10 at 17.28pt

\font\twelvebsy=cmbsy10 at 12pt
\font\tenbsy=cmbsy10
\font\eightbsy=cmbsy8
\font\sevenbsy=cmbsy7
\font\sixbsy=cmbsy6
\font\fivebsy=cmbsy5

\font\tenmsa=msam10

\font\sevenmsa=msam7
\font\fivemsa=msam5
\newfam\msafam
  \textfont\msafam=\tenmsa
  \scriptfont\msafam=\sevenmsa
  \scriptscriptfont\msafam=\fivemsa

\font\tenmsb=msbm10
\font\eightmsb=msbm8
\font\sevenmsb=msbm7
\font\fivemsb=msbm5
\newfam\msbfam
  \textfont\msbfam=\tenmsb
  \scriptfont\msbfam=\sevenmsb
  \scriptscriptfont\msbfam=\fivemsb
\def\Bbb{\fam\msbfam\tenmsb}

\font\tenCal=eusm10
\font\sevenCal=eusm7
\font\fiveCal=eusm5
\newfam\Calfam
  \textfont\Calfam=\tenCal
  \scriptfont\Calfam=\sevenCal
  \scriptscriptfont\Calfam=\fiveCal
\def\Cal{\fam\Calfam\tenCal}

\font\teneuf=eusm10
\font\teneuf=eufm10
\font\seveneuf=eufm7
\font\fiveeuf=eufm5
\newfam\euffam
  \textfont\euffam=\teneuf
  \scriptfont\euffam=\seveneuf
  \scriptscriptfont\euffam=\fiveeuf

\font\seventeenbfit=cmmib10 at 17.28pt

\font\twelvebfit=cmmib10 at 12pt
\font\tenbfit=cmmib10
\font\eightbfit=cmmib8
\font\sevenbfit=cmmib7
\font\sixbfit=cmmib6
\font\fivebfit=cmmib5
\newfam\bfitfam
  \textfont\bfitfam=\tenbfit
  \scriptfont\bfitfam=\sevenbfit
  \scriptscriptfont\bfitfam=\fivebfit


\catcode`\@=11
\def\eightpoint{%
  \textfont0=\eightrm \scriptfont0=\sixrm \scriptscriptfont0=\fiverm
  \def\rm{\fam\z@\eightrm}%
  \textfont1=\eighti \scriptfont1=\sixi \scriptscriptfont1=\fivei
  \def\oldstyle{\fam\@ne\eighti}%
  \textfont2=\eightsy \scriptfont2=\sixsy \scriptscriptfont2=\fivesy
  \textfont\itfam=\eightit
  \def\it{\fam\itfam\eightit}%
  \textfont\slfam=\eightsl
  \def\sl{\fam\slfam\eightsl}%
  \textfont\bffam=\eightbf \scriptfont\bffam=\sixbf
  \scriptscriptfont\bffam=\fivebf
  \def\bf{\fam\bffam\eightbf}%
  \textfont\ttfam=\eighttt
  \def\tt{\fam\ttfam\eighttt}%
  \textfont\msbfam=\eightmsb
  \def\Bbb{\fam\msbfam\eightmsb}%
  \abovedisplayskip=9pt plus 2pt minus 6pt
  \abovedisplayshortskip=0pt plus 2pt
  \belowdisplayskip=9pt plus 2pt minus 6pt
  \belowdisplayshortskip=5pt plus 2pt minus 3pt
  \smallskipamount=2pt plus 1pt minus 1pt
  \medskipamount=4pt plus 2pt minus 1pt
  \bigskipamount=9pt plus 3pt minus 3pt
  \normalbaselineskip=9pt
  \setbox\strutbox=\hbox{\vrule height7pt depth2pt width0pt}%
  \let\bigf@ntpc=\eightrm \let\smallf@ntpc=\sixrm
  \normalbaselines\rm}
\catcode`\@=12

\def\eightpointbf{%
 \textfont0=\eightbf   \scriptfont0=\sixbf   \scriptscriptfont0=\fivebf
 \textfont1=\eightbfit \scriptfont1=\sixbfit \scriptscriptfont1=\fivebfit
 \textfont2=\eightbsy  \scriptfont2=\sixbsy  \scriptscriptfont2=\fivebsy
 \eightbf
 \baselineskip=10pt}

\def\tenpointbf{%
 \textfont0=\tenbf   \scriptfont0=\sevenbf   \scriptscriptfont0=\fivebf
 \textfont1=\tenbfit \scriptfont1=\sevenbfit \scriptscriptfont1=\fivebfit
 \textfont2=\tenbsy  \scriptfont2=\sevenbsy  \scriptscriptfont2=\fivebsy
 \tenbf}

\def\twelvepointbf{%
 \textfont0=\twelvebf   \scriptfont0=\eightbf   \scriptscriptfont0=\sixbf
 \textfont1=\twelvebfit \scriptfont1=\eightbfit \scriptscriptfont1=\sixbfit
 \textfont2=\twelvebsy  \scriptfont2=\eightbsy  \scriptscriptfont2=\sixbsy
 \twelvebf
 \baselineskip=14.4pt}

\def\seventeenpointbf{%
 \textfont0=\seventeenbf  \scriptfont0=\twelvebf  \scriptscriptfont0=\eightbf
 \textfont1=\seventeenbfit\scriptfont1=\twelvebfit\scriptscriptfont1=\eightbfit
 \textfont2=\seventeenbsy \scriptfont2=\twelvebsy \scriptscriptfont2=\eightbsy
 \seventeenbf
 \baselineskip=20.736pt}


\newdimen\srdim \srdim=\hsize
\newdimen\irdim \irdim=\hsize
\def\NOSECTREF#1{\noindent\hbox to \srdim{\null\dotfill ???(#1)}}
\def\SECTREF#1{\noindent\hbox to \srdim{\csname REF\romannumeral#1\endcsname}}
\def\INDREF#1{\noindent\hbox to \irdim{\csname IND\romannumeral#1\endcsname}}
\newlinechar=`\^^J
\def\openauxfile{
  \immediate\openin1\jobname.aux
  \ifeof1
  \message{^^JCAUTION\string: you MUST run TeX a second time^^J}
  \let\sectref=\NOSECTREF \let\indref=\NOSECTREF
  \else
  \input \jobname.aux
  \message{^^JCAUTION\string: if the file has just been modified you may
    have to run TeX twice^^J}
  \let\sectref=\SECTREF \let\indref=\INDREF
  \fi
  \message{to get correct page numbers displayed in Contents or Index
    Tables^^J}
  \immediate\openout1=\jobname.aux
  \let\END=\end \def\end{\immediate\closeout1\END}}

\newbox\titlebox   \setbox\titlebox\hbox{\hfil}
\newbox\sectionbox \setbox\sectionbox\hbox{\hfil}
\def\folio{\ifnum\pageno=1 \hfil \else \ifodd\pageno
           \hfil {\eightpoint\copy\sectionbox\kern8mm\number\pageno}\else
           {\eightpoint\number\pageno\kern8mm\copy\titlebox}\hfil \fi\fi}
\footline={\hfil}
\headline={\folio}

\def\titlerunning#1{\setbox\titlebox\hbox{\eightpoint #1}}
\def\title#1{\noindent\hfil$\smash{\hbox{\seventeenpointbf #1}}$\hfil
             \titlerunning{#1}\medskip}

\newcount\numbersection \numbersection=-1
\def\sectionrunning#1{\setbox\sectionbox\hbox{\eightpoint #1}
  \immediate\write1{\string\def \string\REF
      \romannumeral\numbersection \string{%
      \noexpand#1 \string\dotfill \space \number\pageno \string}}}
\def\section#1{%
  \par\vskip0.666cm\penalty -100
  \vbox{\baselineskip=14.4pt\noindent{{\twelvepointbf #1}}}
  \vskip2pt
  \penalty 500
  \advance\numbersection by 1
  \sectionrunning{#1}}

\def\subsection#1{%
  \par\vskip0.5cm\penalty -100
  \vbox{\noindent{{\tenpointbf #1}}}
  \vskip1pt
  \penalty 500}

\newcount\numberindex \numberindex=0
\def\index#1#2{%
  \advance\numberindex by 1
  \immediate\write1{\string\def \string\IND #1%
     \romannumeral\numberindex \string{%
     \noexpand#2 \string\dotfill \space \string\S \number\numbersection,
     p.\string\ \space\number\pageno \string}}}

\newdimen\itemindent \itemindent=\parindent

\def\item#1{\par\noindent\hangindent\itemindent%
            \rlap{#1}\kern\itemindent\ignorespaces}
\def\itemitem#1{\par\noindent\hangindent2\itemindent%
            \kern\itemindent\rlap{#1}\kern\itemindent\ignorespaces}
\def\itemitemitem#1{\par\noindent\hangindent3\itemindent%
            \kern2\itemindent\rlap{#1}\kern\itemindent\ignorespaces}

\long\def\claim#1|#2\endclaim{\par\vskip 5pt\noindent
{\tenpointbf #1.}\ {\it #2}\par\vskip 5pt}

\def\proof{\noindent{\it Proof}}

\def\today{\ifcase\month\or
January\or February\or March\or April\or May\or June\or July\or August\or
September\or October\or November\or December\fi \space\number\day,
\number\year}

\catcode`\@=11
\newcount\@tempcnta \newcount\@tempcntb
\def\timeofday{{%
\@tempcnta=\time \divide\@tempcnta by 60 \@tempcntb=\@tempcnta
\multiply\@tempcntb by -60 \advance\@tempcntb by \time
\ifnum\@tempcntb > 9 \number\@tempcnta:\number\@tempcntb
  \else\number\@tempcnta:0\number\@tempcntb\fi}}
\catcode`\@=12

\def\bibitem#1&#2&#3&#4&%
{\hangindent=0.8cm\hangafter=1
\noindent\rlap{\hbox{\eightpointbf #1}}\kern0.8cm{\rm #2}{\it #3}{\rm #4.}}


\def\bC{{\Bbb C}}

\def\bQ{{\Bbb Q}}

\def\bZ{{\Bbb Z}}


\def\cO{{\Cal O}}


\def\square{{\hfill \hbox{
\vrule height 1.453ex  width 0.093ex  depth 0ex
\vrule height 1.5ex  width 1.3ex  depth -1.407ex\kern-0.1ex
\vrule height 1.453ex  width 0.093ex  depth 0ex\kern-1.35ex
\vrule height 0.093ex  width 1.3ex  depth 0ex}}}
\def\qed{\kern10pt$\square$}
\def\hexnbr#1{\ifnum#1<10 \number#1\else
 \ifnum#1=10 A\else\ifnum#1=11 B\else\ifnum#1=12 C\else
 \ifnum#1=13 D\else\ifnum#1=14 E\else\ifnum#1=15 F\fi\fi\fi\fi\fi\fi\fi}
\def\msatype{\hexnbr\msafam}
\def\msbtype{\hexnbr\msbfam}
\mathchardef\restriction="3\msatype16   
\mathchardef\boxsquare="3\msatype03
\mathchardef\preccurlyeq="3\msatype34
\mathchardef\compact="3\msatype62
\mathchardef\smallsetminus="2\msbtype72   
\mathchardef\subsetneq="3\msbtype28
\mathchardef\supsetneq="3\msbtype29
\mathchardef\leqslant="3\msatype36   
\mathchardef\geqslant="3\msatype3E   
\mathchardef\stimes="2\msatype02
\mathchardef\ltimes="2\msbtype6E
\mathchardef\rtimes="2\msbtype6F

\def\dbar{\overline\partial}
\def\ddbar{\partial\overline\partial}

\let\ol=\overline

\let\wt=\widetilde

\let\text=\hbox
\def\buildo#1^#2{\mathop{#1}\limits^{#2}}
\def\buildu#1_#2{\mathop{#1}\limits_{#2}}
\def\ort{\mathop{\hbox{\kern1pt\vrule width4.0pt height0.4pt depth0pt
    \vrule width0.4pt height6.0pt depth0pt\kern3.5pt}}}

\def\vlra{\mathrel{\smash-}\joinrel\mathrel{\smash-}\joinrel%
    \kern-2pt\longrightarrow}
\def\srelbar{\vrule width0.6ex height0.65ex depth-0.55ex}
\def\merto{\mathrel{\srelbar\kern1.3pt\srelbar\kern1.3pt\srelbar
    \kern1.3pt\srelbar\kern-1ex\raise0.28ex\hbox{${\scriptscriptstyle>}$}}}

\def\Re{\mathop{\rm Re}\nolimits}
\def\Im{\mathop{\rm Im}\nolimits}

\def\rank{\mathop{\rm rank}\nolimits}

\def\Vol{\mathop{\rm Vol}\nolimits}
\def\Nklt{\mathop{\rm Nklt}\nolimits}
\def\lct{\mathop{\rm lct}\nolimits}


\long\def\InsertFig#1 #2 #3 #4\EndFig{\par
\hbox{\hskip #1mm$\vbox to#2mm{\vfil\special{"
(/home/demailly/psinputs/grlib.ps) run
#3}}#4$}}
\long\def\LabelTeX#1 #2 #3\ELTX{\rlap{\kern#1mm\raise#2mm\hbox{#3}}}


\itemindent = 7mm

\title{Bergman kernels}
\title{and the pseudoeffectivity of}
\title{relative canonical bundles}
\title{-I-}
\titlerunning{Bergman kernels and the pseudoeffectivity of  relative canonical bundles}
\vskip10pt

{\noindent\hangindent0.6cm\hangafter-2
{\bf Bo BERNDTSSON${}^1$ \hfill Mihai P\u AUN${}^2$\kern0.45cm\\
}}

{\noindent\hangindent0.6cm\hangafter-4{\it
\llap{${}^1$}Chalmers University of Technology\hfill
\llap{${}^2$}Universit\'e Henri Poincar\'e\kern0.6cm\break
 Departement of Mathematics\hfill Institut \'Elie Cartan\kern0.6cm\break
   S-41296 G\"oteborg \hfill 54000 Nancy\kern0.6cm\break
}}
\vskip20pt

\noindent{\bf Abstract. \it {The main result of the present article is a (practically optimal)
criterion for the pseudoeffectivity of the twisted relative canonical bundles of surjective projective maps.
Our theorem has several applications in algebraic geometry; to start with, we obtain the natural analytic generalization of some semipositivity results due to E. Viehweg and F. Campana. As a byproduct, we give
 a simple and direct proof of a recent result due to
 C. Hacon--J. McKernan and S. Takayama concerning  
the extension of twisted pluricanonical forms. More applications will
be offered in the sequel of this article.}  } 
\vskip20pt

\section{\S0 Introduction}

In this article our primary goal is to establish some positivity results 
concerning the twisted relative canonical bundle of projective morphisms.

Let $X$ and $Y$ be non-singular projective manifolds, and let 
$p: X\to Y$ be a surjective projective map, whose relative dimension is equal to $n$. 
Consider also a line bundle 
$L$ over $X$, endowed with a -possibly singular- metric $h=e^{-\phi}$, such that the
curvature current  
is semipositive. We denote by ${\cal I}(h)$ the
multiplier ideal sheaf of $h$ (see e.g. [10], [21], [25]). Let $X_y$
be the fiber of $p$ 
over a point $y\in Y$,  such that $y$ is not a critical value of
$p$. We also assume at first   that 
the restriction of the metric $\phi$ to $X_y$  is
not identically 
$-\infty$. 
Under these circumstances, the space of $(n,0)$ forms $L$-valued on $X_y$  
which belong to the multiplier ideal sheaf of the restriction of the metric
$h$ is endowed with a natural $L^2$--metric as follows
$$\Vert u\Vert^2_y:= \int_{X_y}c_nu\wedge\overline u\exp(-\phi)$$
(we use the standard abuse of notation in the relation above).
\noindent Let us consider an orthonormal basis $(u_j^{(y)})$ of the space 
$H^0\bigl(X_y, (K_{X_y}+ L)\otimes 
{\cal I}(h)\bigr)$, endowed with the $L^2$ metric above. 
The relative canonical bundle of $p$ is denoted 
 $K_{X/Y}:=  K_X-p^*K_Y$. 
Recall that the bundles $\displaystyle K_{X_y}$ and $\displaystyle K_{X/Y}\vert X_y$
are isomorphic. Via this identification (which will be detailed in the paragraph 1)
the sections above
can be used to define a metric on the bundle $\displaystyle K_{X/Y}+ L$
restricted to the fiber $X_y$, called  
the {\sl Bergman kernel metric}. This definition immediately extends
also to fibers such that the metric $\phi$ {\it is} identically equal to
$-\infty$ on the fiber. In this case the Bergman kernel vanishes
identically on the fiber, and the Bergman kernel metric is also equal
to $-\infty$ there.

Let $Y^0\subset Y$ be the Zariski open set of points
that are not critical values of $p$ in $Y$, and 
let $X^0\subset X$ be the inverse image of $Y^0$ with respect to $p$.
As $y$ varies in $Y^0$, the above construction defines the 
{\sl relative Bergman kernel metric} on the $\displaystyle K_{X/Y}+ L$ over $X^0$.

Then we have the next result, which gives a 
pseudoeffectivity criterion for the bundle $K_{X/Y}+ L$. 

\claim 0.1 Theorem| Let $p:X\to Y$ be a surjective projective map
between smooth manifolds, and let $(L, h)$
be a holomorphic line bundle endowed with a metric $h$ such that:
\smallskip 
\item {(1)} the curvature current of $(L, h)$ is semi-positive on $X$, i.e. $\Theta_h(L)\geq 0$;
\smallskip
\item {(2)}
$\displaystyle H^0\bigl(X_y, (K_{X_y}+ L)\otimes 
{\cal I}(h)\bigr)\neq 0$
 for some  $y\in Y^0$. 
 
 \noindent Then the relative Bergman kernel metric  of the bundle
$\displaystyle K_{X/Y}+ L\vert X^0$ is not identically $-\infty$. It has
semipositive curvature current 
and extends across $X\setminus X^0$ to a 
metric with semipositive curvature current on all of $X$. 
\endclaim

Several versions of the theorem above were established by the 
first  author in his series of articles on the plurisubharmonic
variation of the Bergman kernels (see [1], [2], [3] and also [22] for
the first results in this direction). 
Let us point out the main
improvements we have got in the present article. 
In the first place, we allow the metric $h$
to be singular. Secondly, the map $p$ is not supposed to be   
a smooth fibration--this will be crucial for the applications, as we will see
in a moment. 
The way we are dealing with the singularities of $p$
is by a careful estimate of the local weight of the 
fiberwise Bergman kernel metric near the singular points: the 
Ohsawa-Takegoshi extension theorem shows that
its local weights are uniformly bounded as we are getting close to the
singular loci of $p$. These local weights are plurisubharmonic 
in $X_0$, so it follows from classical pluripotential theory that they
extend uniquely to  plurisubharmonic functions across  $X\setminus X_0$. This
means that the relative Bergman kernel metric extends uniquely to a
metric with semipositive curvature current across $X\setminus X_0$. 

A similar result has been announced by H. Tsuji in [31] and [31b]. His
idea of proof 
is based on an interesting reduction to the case of a locally trivial
fibration (and then use of the results from [2]), but the somehow sketchy
argument does not seem to be quite complete.
\medskip

One aspect of theorem 0.1 can be seen as a "global"
version of the Ohsawa-Takegoshi theorem: indeed, if $Y$ is just a 
small polydisk, then under the hypothesis (1) and (2) above
the theorem of Ohsawa-Takegoshi shows the existence of 
a section of the adjoint bundle extending the one we have on the
central fiber, thus we get more than the pseudo-effectivity of the adjoint
bundle. 
In the case under consideration, the base $Y$ is compact and simple examples 
show that we cannot expect such a statement, but 
the theorem 0.1 implies that the metric version 
(i.e. replacing effectivity by pseudo-effectivity) still holds.
\vskip 10pt
 
\noindent Assume that $E$ is a pseudoeffective line bundle
on a projective manifold $X$, in the sense that $E$ carries a (possibly singular)
metric with semipositive curvature. Then as a consequence of the 
$L^2$ theory, there exists an 
ample line bundle $A$ on $X$ such that $H^0~(X, mE+ A)\neq 0$, for all $m\gg 0$
(see e.g. [10]).
In this context, as a first application of the theorem 0.1 we have the next
statement.

\claim 0.2 Corollary|Let $p:X\to Y$ be a projective surjective map
between non-singular manifolds $X$ and $Y$. Let  
$(L, h_L)$ be a pseudo-effective line bundle over $X$. Assume that for some
$y\in Y$, not a critical value of $p$, we have $\displaystyle H^0\bigl(
X_y, (K_{X_y}+ L)\otimes {\cal I}(h_L)
\bigr)\neq 0$. Then
for any ample line bundle $A\to Y$ which is positive enough we have 
$H^0\bigl(X, m(K_{X/Y}+ L)
+ p^*A\bigr)\neq 0$, for all $m\geq 0$. 
\endclaim

\noindent Of course, the content of the above statement is that
the positivity we have to add to the bundle $m(K_{X/Y}+ L)$
in order to make it effective comes from $Y$. 
A particular case of this result can be derived from the 
work of F. Campana (see [6]); also, it is consistent with the
semi-positivity results  obtained by T. Fujita, Y. Kawamata, J. Koll\'ar and E. Viehweg
(see [15], [19], [20], [37] as well as the references therein). 

Let us explain in a few words how the corollary is derived 
from  theorem 0.1. 
We have the  decomposition
$$m(K_{X/Y}+ L)
+ p^*A= K_{X}+ (m- 1)(K_{X/Y}+ L)+ L+ p^*(A-K_Y).$$
We denote by $h_B$ the metric obtained in the theorem 1; we use it
to endow the bundle $(m- 1)(K_{X/Y}+ L)+ L+ p^*(A-K_Y)$
with the metric $\displaystyle h_B^{m-1}\otimes h_L\otimes h_{p^*(A-K_Y)}$
where $A\to Y$ is assumed to be positive enough, to compensate for the
possible negativity of the canonical bundle of $Y$ (for the precise positivity
properties of $A$, see the paragraph 3). Now by hypothesis there exists 
$\displaystyle u_y\in H^0\bigl(X_y, (K_{X_y}+ L)\otimes 
{\cal I}(h)\bigr)$ and if we denote by $u_A$ some section of 
$p^*A$, then the section $u_y^{\otimes m}\otimes u_A$
satisfies the integrability properties needed in order to extend it over $X$
by using once again the Ohsawa-Takegoshi theorem. 
The fact that the metric $h_B$ is explicitely given over 
a Zariski dense open set is  crucial here.

\vskip 15pt
In Corollary 0.2 we prove the existence of global sections to $m(K_{X/Y}+L)
+A$, assuming the existence of sections to $(K_{X/Y}+L)\otimes
{\cal I}(h)$ over at least 
some fiber. It is natural to ask what happens if we only assume from
the start existence of certain sections over fibers to multiples of
$K_{X/Y}+L$. In  section 4 we prove a version of Theorem 0.1
for the {\it Narasimhan-Simha } metric, see [26]
, which is an analog of
the Bergman kernel metric for multiples of the canonical bundle,
recently reviewed by Tsuji [32], [32a] and [32b]. This
generalizes a classical result of Kawamata [18] for the nontwisted
case and will be developed further in the sequel to this paper [4]. 

The proofs in section 4 are still based on the same positivity results
for direct image bundles from [3] as the proof of Theorem 0.1, but require an
additional twist. In particular, we need  (weaker) versions of the
results from [3] for nonsmooth metrics. These are discussed in section 3 -
hopefully they are also of some independent interest.

\vskip 10pt
To introduce the next application, let us recall the notion of 
{\sl restricted volume} (see [5], [14]). Let $E\to X$ be a holomorphic line bundle.
If $V\subset X$ is an irreducible
$d$-dimensional sub-manifold, let us denote by
$$H^0(X\vert V, mE):= \Im (H^0(X, mE)\to H^0(V, mE_V)),$$
and let $h^0(X\vert V,mE)$ be the dimension of this space.
Then the restricted volume of $E$ to $V$ is
$$\Vol_V(X,L):= \lim\sup {{d!}\over {m^d}} h^0(X\vert V, mL))$$

The definition of a {\sl maximal center} is slightly more involved,
and it will not be recalled here (see e.g. [16], [30] and also the paragraph 6
of this article). Let us just mention that 
given an effective $\bQ$--divisor $D= \sum_j\nu_jZ_j$, the maximal
centers of the pair  
$(X, D)$ are the higher codimensional analog of the 
irreducible components $Z_j$ above
such that $\nu_j= 1$.

The corollary above can be used 
to provide a rather  simple  proof
of the next statement due to
S. Takayama ([30], see also [16] and [32]), which is crucial in 
the investigation of the properties of the 
pluricanonical series. 

\claim 0.3 Theorem {\rm ([29], [30], [16])}| Let $X$ be a non-singular
projective manifold, 
such that $K_X$ is pseudo-effective. Let
$L$  be a line bundle which can be written as $L= A+ D$ where 
$A$ is $\bQ$-ample and $D$ is a $\bQ$-effective divisor. Assume that
$V$ is an irreducible maximal center of 
$(X, D)$ such that there exists a section of some multiple of $K_X+L$ 
vanishing on some ample divisor and
which is not identically zero when restricted to $V$. Then 
$$\Vol_V(X,K_X+ L)\geq \Vol(V).$$
\endclaim

Even if our proof goes along the same
lines as the previous ones, the exposition is substantially simplified by
the use of the corollary 0.2,which allows us to bypass the use of  Kawamata's
subadjunction theorem [19].

There are 2 major ingredients
needed in the proof of 0.3. The first one is the technique invented by
Y.-T. Siu to  
prove the invariance of plurigenera, see [33], [34]. The version which will be
used here is due to S. Takayama, but we will offer a simpler proof, in the same
spirit as in [28]. The second one is simply the corollary 0.2.

Let us explain  vaguely how the two techniques
combine to prove the theorem. First, if
the maximal center $V$ is a hypersurface, then the
theorem 0.3 is a consequence of the invariance of plurigenera techniques.
If $V$ has higher codimension, one is lead  to reduce to the
divisor case by using a modification $\mu: Y\to X$.
We consider the restriction of the map $\mu$
to a well-chosen exceptional divisor $S$ which maps onto
$V$. We know how to
extend  twisted pluricanonical sections on $S$.
The crucial result which allows us to inject the space of
pluricanonical sections on $V$ to a space of twisted pluricanonical
sections on $S$ is precisely Corollary 0.2 - the injection map is 
given by multiplication with a section of a twisted relative canonical
bundle $K_{S/V}+B$, whose existence is guaranteed by Corollary 0.2.

\vskip 10pt
\claim Acknowledgements|The second named author  would like to
thank F. Campana, Y. Kawamata, S. Takayama and E. Viehweg
for very interesting discussions concerning various aspects of this
work; he would equally like to express his gratitude to his friends
and colleagues 
at the Institute \'Elie Cartan (Nancy) for the nice atmosphere during
the ``groupe de travail" 
where topics like the extension of pluricanonical forms were
extensively analyzed.  We also would like to thank the referees for
detailed and constructive criticism of the original manuscript.

\endclaim

\vskip 30pt
\section{\S1 Relative Bergman kernels}  
\vskip 5pt
\noindent Let 
 $p:X \mapsto Y$ be a surjective holomorphic map with 
compact $n$-dimensional fibers, onto a complex manifold $Y$. At first we will
 also assume that the map $p$ defines a smooth fibration, i e that the
 differential $dp$ is surjective at each point, so that the fibers
$$
X_y:=p^{-1}(y)
$$
are smooth manifolds. Let  $(L,\phi)$ be a Hermitian holomorphic line
bundle over the total space $X$. For each $y$ in $Y$ we let
$$
E_y=H^0(X_y,L|_{X_y}+K_{X_y})
$$
be the space of global holomorphic $L$-valued $(n,0)$-forms on
$X_y$. The Hermitian metric $\phi$ on $L$ induces a Hilbert norm $H_y$
on $E_y$ by
$$
\Vert u\Vert_y^2 =\int_{X_y} c_n u\wedge\overline u e^{-\phi}. \leqno(1.1)
$$

Near a point in $Y$ we can choose local coordinates $t=(t_1, ...t_m)$,
which define a local $(m,0)$-form, $dt=dt_1\wedge ...\wedge dt_m$, that
trivializes the canonical bundle of $Y$, $K_Y$. A choice of such form
gives a natural map from forms, $u$, on the fibers $X_y$
to sections, $\tilde u$, of $K_X$ over $X_y$ by
$$
\tilde u= u\wedge p^*(dt)=u\wedge dt,
$$
where in the last equality we abuse notation slightly by thinking of
$t$ as a function on $X$. When dealing with forms on $X$ we will
mostly use
the notation $dt$, but we will recurr to  $p^*(dt)$ for emphasis occasionally.

Conversely, given a  local section $\tilde u$ of $K_X$ we can write
$\tilde u=u\wedge dt$ 
locally. The restriction of $u$ to fibers is then uniquely defined and
thus defines a  form on $K_{X_y}$. We will call $u$ the trace of
$\tilde u$ on $K_{X_y}$. Altogether the correspondence between $u$
and $\tilde u$ gives us an identification between forms on a fiber and
restrictions of sections to $K_X$ to the fiber. 
This identification clearly depends on the choice of $dt$, but the map
$ (u,dt)\mapsto u\wedge dt$
is  an invariantly defined isomorphism
$$
K_{X_y}+p^*(K_Y)|_{X_y}\mapsto  K_X|_{X_y},
$$
for any $y$. Hence we also get an isomorphism
$$
K_{X_y}\mapsto  (K_X-p^*(K_Y))|_{X_y} =K_{X/Y}|_{X_y}, 
$$
where $K_{X/Y}=(K_X-p^*(K_Y))$ is the {\it relative canonical bundle}
of $p$. If $t$ is any choice of local coordinates on $Y$, then the
isomorphism is given by
$$
u\mapsto u\wedge p^*(dt)/dt.
$$
The Hilbert norm (1.1) therefore defines a norm on sections of $ L+
K_{X/Y}$ over $X_y$.

\medskip

\noindent For any $y$ in $Y$ we now define the Bergman kernel for this norm by
$$
B_y =\sum u_j\otimes\overline {u_j},
$$
for any choice of orthonormal basis $(u_j)$ of the space of sections of
$(L+K_{X/Y})|_{X_y}$. It is well known 
that this definition is independent of the choice of orthonormal
basis. In the manifold case, that we are dealing with here, the
Bergman kernel is not a function, but $\log B_y$ defines a metric on
 $(L+K_{X/Y})|_{X_y}$ in the sense that
$$
u\otimes\overline u/ B_y
$$
is a well defined function if $u$ is a section of  $L+K_{X/Y}$ over
$X_y$. Alternately, we could choose an arbitrary  smooth metric, $\chi$,
on  $(L+K_{X/Y})|_{X_y}$. Then $B_y e^{-\chi}=:e^{\xi}$ is a well defined global
function on $X_y$. Hence
$$
\log B_y=\xi+\chi
$$
is also a metric on  $(L+K_{X/Y})|_{X_y}$.

As $y$ varies, we get a metric on the line bundle  $L+K_{X/Y}$
over all of $X$.
We will refer to this metric as the relative Bergman kernel metric and
write $B= e^{\psi}$, with the understanding that $\psi$ is a metric
which is represented by different local functions, given different
local trivializations. 

In the classical case of domains in $\bC^n$ there is a well known extremal
characterization of the Bergman kernel as
$$
B(x)=\sup_{\|h\|\leq 1}|h(x)|^2.
$$
We will use repeatedly a variant of this in our setting. Choose local
coordinates $x$ near a point in a fiber $X_y$. Then $dx/dt$ is a local
frame for $K_{X/Y}$. Choose also some
local trivialization, $e$, of $L$ near the same point. With respect to these
local frames $\psi$ is represented by a function $\psi'$, and any
section $\tilde u/dt=u\wedge p^*(dt)/dt$ to  $L+K_{X/Y}$ over $X_y$
can be written 
$$
\tilde u/dt=u'(x)dx/dt\otimes e.
$$
Then 
$$
e^{\psi'(x)}=\sup_{\|u\|_y\leq 1} |u'(x)|^2.
$$
The proof of this is exactly as in the classical case. It is of
paramount importance that, since our local trivializations are given by
$dx/dt$, the coefficient that enters in the extremal characterization
is the coefficient of $\tilde u$ (and not any coefficient of $u$). It
is the presence of the factor $p^*(dt)$ - which vanishes on the
singular locus of $p$ - in $\tilde u$ that will allow
us to control 
the Bergman kernel metric near the singular locus.

Our first  result, Theorem 0.1, says that if $\phi$ is a, possibly
singular, metric of nonnegative curvature on $L$, then $\psi$ is a
metric of nonnegative curvature of  $L+K_{X/Y}$ - unless $B$ is
identically 0. More precisely, $\psi$ is semipositive over the Zariski
open set where $p$ is a smooth fibration, and extends in a unique way
to a semipositive metric over the singular fibers.

To start explaining the proof of this, let us first assume that $\phi$
is smooth, and that $p$ defines a smooth fibration, i e that $dp$ is
surjective everywhere. Theorem 0.1 is then an easy consequence of
Theorem 1.2 in [3]. This theorem  says that the Hilbert
norms $\| \cdot\|_y$ on $E_y$ define a semipositively curved Hermitian
metric on a vector bundle $E$ over $Y$, with fibers $E_y$. The complex
structure on $E$ is such that a local section of $E$
$$
y\mapsto u_y
$$
is holomorphic if and only if the associated sections
$$
y\mapsto \tilde u_y =u_y\wedge dt
$$
to $K_X$ are holomorphic, and hence define holomorphic sections
$\tilde u/dt$ to $K_{X/Y}$. ( The bundle $E$ is the direct image
$p_*(L+K_{X/Y})$ of $L+K_{X/Y}$ under the map $p$.) The precise
statement of Theorem 1.2 in [3] is that $E$ with the $L^2$- Hermitian
metric is (weakly) Nakano positive, but the weaker property of
Griffiths positivity is all we use here.

Choose local coordinates $(t,z)$ near a point in $X$, with respect
to which the map $p$ is the trivial fibration $(t,z)\mapsto t$. These
local coordinates give us local frames, $dt$ and $dt\wedge dz$ for
$p^*(K_Y)$ and $K_X$ respectively, and hence the local frame $dt\wedge
dz/ dt$ for $K_{X/Y}$. With respect to this frame, $\psi$ is given by
a function $\psi'$, and we need to prove that this function is
plurisubharmonic. This means that
$$
\psi'(t,h(t))
$$
is plurisubharmonic of $t$ if $(t,h(t))$ is any local holomorphic map.

We now use the same local coordinates, and the map $h$, to define a local
holomorphic section of the dual bundle of $E$, $E^*$. Let $u_y$ be an
element in $E_y$, and let
$$
u_y\wedge p^*(dt)/dt
$$
be the associated  element in $K_{X/Y}+L$. With respect to the
trivializations
$$
u_y\wedge p^*(dt)/dt=u'(t,z)(dz\wedge dt)/dt\otimes e.
$$
We put
$$
\xi_y(u_y)= u'(t,h(t)).
$$
If $u_y$ depends holomorphically on $y$, then $u'$ is a holomorphic
function of $t$ and $z$, so it follows that $y\mapsto \xi_y$ is a
holomorphic section of $E^*$. By the extremal characterization of
Bergman kernels
$$
\|\xi_y\|^2=\sup_{\|u\|\leq 1} |u'(t,h(t))|^2= e^{\psi'(t,h(t))}.
$$
Since the logarithm of the norm of any section of the dual of a
Griffiths positive bundle is plurisubharmonic, it follows that 
$\psi'$ is plurisubharmonic, which is what we wanted to prove.

This proves Theorem 0.1 in the smooth case, i e when both the metric
$\phi$ and the fibration $p$ are smooth. Let us now still assume that
$\phi$ is smooth but relax the assumption on $p$, so that $p$ is
a general surjective holomorphic map. Then the degeneracy locus of $p$
is an analytic subvariety $W$ of $X$ and $p(W)$ is an analytic
subvariety of $Y$. By Sard's theorem $p(W)$ has zero measure in $Y$,
so it is a proper subvariety, $W'$ of $Y$. Outside of $p^{-1}(W')$,
$p$ is a smooth fibration, and by the previous discussion, $\psi$ has
semipositive curvature there. We want to prove that $\psi$ extends to
a semipositive metric across  $p^{-1}(W')$. Since  $p^{-1}(W')$ is
pluripolar, it suffices for this to prove that $\psi$ stays locally
bounded from above. The next lemma is the crucial step.

\claim 1.1 Lemma|Let $D$ be a polydisk in $\bC^{n+m}$, and let $\phi$
be a plurisubharmonic function in $D$. Let
$$
p: D\mapsto V
$$
be a holomorphic map from $D$ to a bounded open set in $\bC^m$.
Assume 0 is a regular value of $p$ and let  $D_0:=p^{-1}(0)$. Let $u$
be a holomorphic $(n,0)$-form on $D_0$. Let $\tilde 
u=u\wedge p^*(dt)$ where $t=(t_1, ...t_m)$ are standard coordinates on
$\bC^m$ and write
$$
\tilde u=u' dz,
$$
where $z=(z_1, ... z_{n+m})$ are standard coordinates on
$\bC^{n+m}$. Then there is a constant $C_K$ such that  for any
compact subset $K$ of $D$ and   $z$ in $D_0\cap K$,  
$$
|u'(z)|^2\leq C_K \int_{D_0} [u,u]e^{-\phi}.
$$

\endclaim

\proof.  By the Ohsawa-Takegoshi extension theorem there
is a holomorphic $(n+m,0)$-form $\tilde U$ in $D$ such that
$$
\tilde U=\tilde u
$$
on $D_0$ and 
$$
\int_D |\tilde U|^2 e^{-\phi}\leq C\int_{D_0} [u,u] e^{-\phi},
$$
where the constant $C$ only depends on the sup-norm of $p$. (See
section 5, Theorem 5.1, and the comments immediately after that
theorem for an appropriate version of the Ohsawa-Takegoshi extension
theorem.) Since
$\phi$ is locally bounded from above it follows that
$$
\int_K |\tilde U|^2 \leq C'\int_{D_0} [u,u] e^{-\phi}.
$$
The lemma now follows from the meanvalue inequality applied to $\tilde U$.
\hfill \qed
\vskip 7pt
 \claim 1.2 Remark|{\rm The main point in the lemma is that the constant $C_K$
 does not blow 
up as we approach a singular fiber - the Ohsawa-Takegoshi theorem 
implies that the fiberwise Bergman kernel  automatically has the right
normalization to extend over the 
singular fibers of $p$.}
\endclaim
If we choose local coordinates and local trivializations near a point
on a singular fiber of $p$ it follows immediately from the lemma,
together with the extremal characterization of Bergman kernels, that
the relative Bergman kernel stays bounded near a singular fiber. Hence
the metric $\psi$ extends (uniquely!) to a semipositive metric on
$L+K_{X/Y}$. (Notice that the Bergman kernel  is not identically
equal to 0 since by assumption $K_{X/Y}+L$ has some section over some
fiber in $X_0$.)
\claim 1.3 Remark|{\rm 
Two examples may serve to illustrate what happens near the singular
locus. Let us first look at the case when $X$ is the blow-up of $Y$ at
a point, and, say, $L$ is trivial. Then
$$
s:= p^*(dt)/dt
$$
defines a global holomorphic section of the relative canononical
bundle $K_{X/Y}$ ($s$ does not depend on the choice of local frame
$dt$). Since $s$ vanishes to order $(n-1)$ precisely on the 
exceptional divisor $D$, we see that
$$
K_{X/Y}=(n-1)(D),
$$
and that $h_D:=|s|^2$ defines a singular pseudoeffective metric on
$K_{X/Y}$. We claim that $h_D$ equals the relative Bergman kernel
metric in this case. To see this, note that outside of the singular
locus the fibers $X_y$ are just points, and elements, $u$, in $E_y$
are just numbers. Then 
$$
\tilde u/dt=u\wedge p^*(dt)/dt=u' dx/dt,
$$
and to compute the Bergman kernel of this zero-dimensional fiber we
need to take the supremum of $|u'|^2$ over all $u$ of modulus 1. Clearly
this equals $|s'|^2$ where 
$$
s=s'dx/dt,
$$
proving our claim. A very similar argument applies, and the same
conclusion holds, in the more general
situation of a map $p$ with a finite number of preimages over a generic
point in the base.

 The second example illustrates what happens when the dimension of the
 fiber does
 not jump at the singular locus. The perhaps simplest such example is
 when  $n=m=1$ and $p(z_1,z_2)= z_1^2
 -z_2^3$. Then 
the fiber $D_0=\{ z_1^2 =z_2^3\}$ can be parametrized by
$$
\lambda\mapsto (\lambda^3,\lambda^2),
$$
and any holomorphic form $u$ on $D_0$ can be written
$$
u=h(\lambda) d\lambda,
$$
with $h$ holomorphic off 0. Then $u$ has finite $L^2$-norm on $D_0$ if
and only if $h$ extends holomorphically across 0. On the other hand
$$
dp= 2z_1dz_1 - 3z_2^2dz_2,
$$
and
$$
d\lambda = (1/z_2)dz_1 -(z_1/z_2^2)dz_2,
$$
so we get
$$
\tilde u=u\wedge dp =-\lambda^2 h(\lambda)dz_1\wedge dz_2.
$$
Therfore $\tilde u'$ is finite at the origin if and only if $h$ has
a pole of order at most 2 at 0. Since by the extremal
characterization of Bergman kernels, the Bergman kernel at the origin
equals
$$
\sup |\tilde u'(0)|^2/\|u\|^2,
$$
we see that Bergman kernels tend to be small - not big - on singular fibers.}
\endclaim

It now remains to relax the condition that $\phi$ be smooth to obtain
the full proof of Theorem 0.1. We may again assume that $p$ is a
smooth fibration, since the general case then follows from Lemma 1.1
by the same argument as before. We then need to approximate $\phi$ by
a sequence of smooth semipositive metrics. It is here that we
use the assumption that $p:X\mapsto Y$ is projective, which permits us to
reduce the problem to the Stein case. This is the object of the next section.

\vskip 30pt

\section{\S2 Stein and projective fibrations}
\vskip 5pt
We will now discuss the proof of Theorem 0.1 in the general case.
First, since $X$ is assumed to be projective, we 
can remove a divisor and get a Stein submanifold. On this Stein
manifold the line bundle $L$ has a holomorphic section. Removing the
zero divisor of this section we get a smaller
Stein submanifold over which $L$ is trivial. We will prove 
that the relative Bergman kernel metric for this Stein submanifold has
nonnegative curvature. But, since divisors are removable for $L^2-$
holomorphic functions, the relative Bergman kernel metrics for $X$ and for the
Stein submanifold are identical, so Theorem 0.1 follows.

We first consider a smooth Stein fibration over a domain $W$
in $\bC^m$. By this we mean that we are given a Stein manifold
$D$ of dimension $n+m$, together with a holomorphic surjective map $p$
from $D$ to $W$ whose derivative $dp$ is everywhere of maximal
rank. Then the fibers $D_t=p^{-1}(t)$ are $n$-dimensional Stein
manifolds. Let $\phi$ be a plurisubharmonic function on the total
space $D$. 

We let $E_t$ be the space of global holomorphic $(n,0)$-forms $u$ on $D_t$
that lie in $L^2$ in the sense that
$$
\|u\|^2_t = \int_{D_t} c_n u\wedge\bar u e^{-\phi} <\infty.
$$
This is now, in general, an infinite dimensional space, but we can
still define its Bergman kernel
$$
B_t =\sum u_j\otimes\bar u_j,
$$
and a metric $B=e^\psi$ on $K_{D/W}=K_D$ as in the previous section.
If $u$ lies in $E_t$ we write $\tilde u= u\wedge dt$ and 
$$
\tilde u=u' dz,
$$
given local coordinates on $D$. We then have, by the extremal
characterization of Bergman kernels, that
$$
e^{\psi(z)}=\sup|u'(z)|^2,
$$
with respect to the induced local trivializations, where the supremum
is taken over all forms $u$ of norm at most 1.

\claim 2.1 Theorem|In the above situation, the relative
Bergman kernel metric has 
nonnegative curvature, i e 
$$
i\ddbar\psi\geq 0,
$$
if the Bergman kernel is not identically 0. 
\endclaim

The strategy of the proof is as follows. We first prove the theorem in
a 'good'  situation, where the weight function $\phi$ is smooth and
the Stein manifold is of very special type.  It is then clear that
the relative Bergman kernel 
is continuous, so to prove that it is plurisubharmonic it suffices to
verify subharmonicity on complex lines. This means in particular that
we may assume the base is one dimensional. (This is because the
complex structure on our bundle $E$ restricted to a complex line in
the base coincides with the structure we get if we first restrict the
fibration $X\mapsto Y$ to the line and then construct our vector
bundle from there.) Finally, the general case
is reduced to the 'good' case. 

Since $D$ is Stein it has a smooth strictly plurisubharmonic
exhaustion, $\rho$. Let $D^c=\{\rho<c\}$. By Sard's lemma, $D^c$ is
smoothly bounded and strictly pseudoconvex for  all choices of
$c$ outside a closed set of measure 0.

Let $c$ be one such noncritical value and let $C_c$ be the set of all
$t$ in the disk such that $c$ is critical for the restriction
of $\rho$ to $D^c_t:=D^c\cap p^{-1}(t)$.

\claim 2.2 Lemma| $C_c$ is a closed subset of $\Delta$ of zero measure. 
\endclaim

\proof.  Let $t=\tau+i\sigma$. Note that a point $x$ is critical for the
restriction of $\rho$ to $D_t$ if and only if $d\rho$, $d\tau$ and
$d\sigma$ are linearily dependent at $x$ (here we view $t$ as a
function on $X$). If $x$ is not critical for $\rho$ on $X$, $d\rho$ is
not 0 at $x$, so $d\rho|x$ must be a linear combination of $d\tau$ and
$d\sigma$. This in turn means that $d\tau$ and $d\sigma$ are linearily
dependent on the kernel of $d\rho$, so that the restriction of $p$ to
$\partial D^{\rho(x)}$ has a critical point at $x$. The lemma therefore
  follows from Sard's lemma applied to the restriction of $p$ to $
  \partial D^{\rho(x)}$.
     \qed

Take $t_0$ in $\Delta$.  We
shall prove that the fiberwise Bergman kernel for $D^c$ is
logplurisubharmonic in  
$D^c\cap p^{-1}(V)$, if $V$ is a sufficiently small neighbourhood of
$t_0$. 

\claim 2.3 Lemma| There is a holomorphic tangent vector field, $f$, to
$D$ such that
$$
dp(f)=\partial/\partial t.
$$
\endclaim

\proof. 
   Such a field is easily found locally by choosing local
coordinates with respect to which the fibration is trivial, and a
nonholomorphic field $f'$ is
then obtained globally by a partition of unity. Let $w=\dbar f'$. Then
$w$ is a $\dbar$-closed $(0,1)$-form with values in the subbundle of
the holomorphic tangent bundle of $D$ consisting vectors tangent to
fibers ( i e the kernel of $dp$). Since $D$ is Stein, $w=\dbar g$,
where $g$ is a field with values in the same subbundle. Replacing $f'$
by $f=f'-g$ we get a holomorphic tangent vector field on $D$ satisfying
the claim in the lemma.
\qed

Let $\Omega:= D^{c+1/2}_{t_0}$, and consider the flow of our
field $f$ from the lemma above. The flow after time $t-t_0$ maps $\Omega$
to a domain in $D_t$, which contains $D_t^c$ and is contained in
$D^{c+1}_t$ if $t$ lies in $V$, a small enough neighbourhood of
$t_0$. Let $\tilde \Omega$ be the image of $\Omega$ under the
flow for $t$ in $V$. 

We  now follow the arguments from [2]. Since there is a fiberpreserving
biholomorphism between $\Omega\times V$ and $\tilde\Omega$, it first
follows that the Bergman kernel for $\tilde\Omega$ is
logplurisubharmonic, by the argument from [2], section 2. (In [2], this
is proved for the product of domains in Euclidean space, but the same
arguments apply for products of a domain with a manifold, provided we
consider $(n,0)$-forms instead of functions.) We then
replace our weight function $\phi$ by $\phi_j$, where $\phi_j$ form a
sequence of plurisubharmonic weights, all equal to $\phi$ in $D^c$,
and tending to infinity outside the closure of $D^c$.

To prove the logplurisubharmonicity of the Bergman kernel in $D^c\cap
p^{-1}(V)$ we shall use the following version of Lemma 3.1 from [2].

\claim 2.4 Lemma| Let $G_0$ and $G_1$ be Stein manifolds,
with $G_0$ compactly included in $G_1$. Let $\phi_j$ be a
sequence of continuous weight functions in $G_1$ such that
$\phi_j=\phi$ in $G_0$ and $\phi_j$ increases to infinity almost
everywhere in $G_1\setminus G_0$. Let $B_j$ be the
Bergman kernel for the space $A_j$ of holomorphic $(n,0)$-forms $u$ in
$G_1$ with Hilbert norm
$$
c_n\int_{G_1} u\wedge \bar u \, e^{-\phi_j}.
$$
Then for any $z$ in $G_0$,  $B_j(z,z)$ increases to $B(z,z)$, the
Bergman kernel for the closure of the space $A_1$ in the space $A$ of
holomorphic $(n,0)$-forms $u$ in $G_0$  with Hilbert norm
$$
c_n\int_{G_0} u\wedge \bar u \, e^{-\phi}.
$$
\endclaim

\proof. From the extremal characterization of Bergman kernels it is
clear that $B_j$ increases, and it is also clear that $B_j\leq B$.  To
prove the opposite inequality, we take a form $u$ in $A_1$ with norm at
most 1 in $A$, which almost realizes the supremum in the extremal
characterization. By monotone convergence, the norm of $u$ in $A_j$
tends to the norm of $u$ in $A$ as $j$ tends to infinity. Hence $\lim
B_j(z,z)\geq B(z,z)$.
\qed 

We now apply this Lemma 2.4 to the domains $D^{c}_t$ inside of
$\tilde\Omega$, with weight functions $\phi_j=\phi +\psi_j$, with
$\psi_j$ being a sequence of plurisubharmonic functions all equal to 0
in   $D^{c}_t$ and increasing to infinity in $\tilde\Omega$ outside 
the closure of $D^{c}_t$  (suitable convex functions of $\rho$ will do
for $\psi_j$). We claim that if $t$ is taken outside of the
exceptional set $C_c$ of Lemma 2.2, then the any form on $D^c_t$ can
be approximated by forms on $\tilde\Omega_t$, so that the Bergman
kernel $B$ of Lemma 2.4 equals the Bergman kernel on $D^c_t$.

This claim follows from  the next lemma, which is a  variant of
 Lemma 3.4 from [2]. The version  here is however substantially simpler
 since we have taken care to have all the domains involved smoothly
 bounded, by avoiding critical values of $\rho$. This permits us to
 bypass a nontrivial result of Bruna and Burgues that was used in [2].

\claim 2.5 Lemma| Let $\rho$ be a smooth plurisubharmonic exhaustion
function in an $n$-dimensional Stein manifold $G$, and put $G_c=\{\rho
<c\}$. Let $u$ be a holomorphic form of bidegree $(n,0)$ in $G_0$, such
that 
$$
c_n\int_{G_0} u\wedge\bar u \, e^{-\phi}<\infty,
$$
where $\phi$ is a plurisubharmonic weight function on $G$.
Assume that 0 is not a critical value of $\rho$, and take $c>0$. Then
$u$ can be approximated in $L^2$-norm on $G_0$ by a sequence of
$(n,0)$-forms that extend holomorphically to $G_c$. 
\endclaim

\proof.  Since, for any small $\epsilon>0$,  $G_\epsilon$ is Runge with respect
to $G_c$, it suffices to approximate $u$ by forms holomorphic in some
small neighbourhood of $\bar G_0$. For this, we first approximate $u$
by $v$, smooth in a neighbourhood of $\bar G_0$, with $\bar\partial v$
small, by decomposing $u=\sum u_j$, where $u_0$ has compact support in
$G_0$, and each of the other $u_j$ is supported in a small
neighbourhood of a boundary point. We may assume that these small
neighbourhoods of boundary points are starshaped, and then approximate
the corresponding $u_j$ by  $v_j$ obtained from a small dilation (this is
where we use that 0 is noncritical).  Putting $v_0=u_0$ and $v=\sum
v_j$ it follows that $v$ approximates $u$ in $L^2$ and that
$\bar\partial v$ is small in a neighbourhood of $G_0$. To obtain a
holomorphic approximant we finally solve a $\bar\partial$-equation in
$G_\epsilon$. 
\qed

Consider now the sequence of Bergman kernels on $\tilde\Omega_t$ for
the weights $\phi_j$ . By Lemma 2.4 they form an increasing sequence
that on each $D^c_t$ converges to  the Bergman kernel on $D^c_t$, with
weight $\phi$ for
the space of holomorphic forms that can be approximated by forms
holomorphic on $\tilde\Omega_t$. Since the Bergman kernels on
$\tilde\Omega_t$ for 
the weights $\phi_j$ are nonnegatively curved,  the 
upper semicontinuous reularisation of the limit is also a
nonnegatively curved metric (on $K_D$), that we denote
$e^{\Psi_c}$. On the other hand, when $t$ is outside of $C_c$, Lemma
2.5 says that any holomorphic form on $D^c_t$ can be approximated by
forms holomorphic in $\tilde\Omega_t$, so outside of $p^{-1}(C_c)$ the
limit is just the fiberwise Bergman kernel on $D^c_t$. It is not hard
to check (by a normal family argument) that the fiberwise Bergman
kernels are already uppersemicontinuous outside of   $p^{-1}(C_c)$, so
$e^{\Psi_c}$ is just the fiberwise Bergman kernel on $D^c_t$ there, in
particular almost everywhere.

Summing up, we have proved that, over $V$, the relative Bergman kernel
metric on
$D^c$
is nonnegatively curved outside of $p^{-1}(C_c)$ and extends
(uniquely) to a 
nonnegatively curved metric, $\Psi_c$ across $p^{-1}(C_c)$. Since $V$
is a small neighbourhood of an arbitrary point in the disk, the same
thing holds on all of $D_c$. Let
$e^{\psi_c}$ be the fiberwise Bergman kernel for all $t$ in
$\Delta$. A normal family argument shows that
$$
\psi^c\geq \Psi^c,
$$
and, as we know, equality holds outside $C_c$. Moreover we claim that
$$
\psi^{c+\epsilon}\leq \Psi^c,
$$
for $\epsilon>0$. This follows since if $t$ lies in $\Delta$ and
$\epsilon'<\epsilon$ is noncritical for $\rho$ on the fiber over $t$, 
then
$$
 \psi^{c+\epsilon}\leq \psi^{c+\epsilon'}=
\Psi^{c+\epsilon'}\leq\Psi^c.
$$
Both $\psi^c$ and $\Psi^c$ decrease with
$c$, so they must have the same limit as $c$ tends to infinity.
Since $\Psi_c$ is nonnegatively curved for all $c$ this limit is
nonnegatively curved.

However, for each $t$,   the limit of
$e^{\psi_c}$ as $c$ tends to 
infinity on the fiber over $t$,
equals the Bergman kernel for $D_t$. (To see this, note first that the
limit must be larger than the Bergman kernel for $D_t$ since the
Bergman kernel is a decreasing function of the domain. The opposite
inequality follows from the extremal characterization of Bergman
kernels, since the extremals for each domain $D^c_t$ have a
subsequence that converges weakly in $L^2$ to a form of norm at most 1.)

This proves Theorem 2.1 when $\phi$ is smooth. Finally, to remove the
assumption that $\phi$ be smooth we  
again consider regular values $c$ of the exhaustion function. A
general plurisubharmonic $\phi$ can now be written as a decreasing
limit of smooth strictly plurisubharmonic $\phi_\nu$ in a neighbourhood 
of the closure of $D^c$. Applying Theorem 2.1 to $D^c$ with the weight
$\phi_\nu$, we see that the corresponding Bergman kernels are
plurisubharmonic and decrease to the Bergman kernel for the weight
$\phi$. Hence the limit is again plurisubharmonic, and letting $c$
tend to infinity we see that $\psi$ is plurisubharmonic or
identically equal to $-\infty$ in $D$.  Note that there is no need
to take upper semicontinuous regularizations - it follows from the
argument above that the fiberwise Bergman kernels as they are give us
an uppersemicontinuous metric. This is of course under the
assumption 
in Theorem 2.1, that the fibration $p$ be smooth. We loose the
explicit character of our metric when we extend across the
singularities of the fibration using Lemma 1.1.

Notice also that it may well happen that the limit of the sequence of
Bergman kernel on $D^c_t$ as $c$ tends to infinity is equal to 0. This
happens precisely when there are no nontrivial $L^2$ holomorphic forms
on the fiber $D_t$. If on some fiber there is a nontrivial $L^2$
holomorphic form, then the Bergman kernel on that fiber is not
identically 0, so the relative Bergman kernel metric is a
nonnegatively curved (singular) metric. 
\hfill \qed

Thus, by the discussion at the beginning of this section, the 
proof of the theorem 0.1 is also finished

\section{\S 3 Positivity of the direct image bundle.}

As noted in section 1,  positivity of the relative Bergman kernel metric
(i. e. plurisubharmonicity of the Bergman kernel) follows in the
smooth case from the (Griffiths) semipositivity of the direct image bundle 
$$
E=p_*(K_{X/Y}+L).
$$
A priori, the positivity of $E$ is however a stronger statement: It
implies plurisubharmonicity of {\it any} function $\log|\xi|^2$, for $\xi$ any
holomorphic section of the dual bundle $E^*$, whereas
plurisubharmonicity of Bergman kernels corresponds to the special case
of point evaluations. In this section we shall
prove a variant of this stronger statement when the metric on $L$ is
not necessarily smooth, and also show that it corresponds to a
stronger positivity condition on the Bergman kernel.

To start with we shall discuss the perhaps not completely standard
notion of a singular metric on a vector bundle (see also de Cataldo,
[7]). By this we basically mean a measurable map from the base, $Y$, to the
space of nonnegative hermitian forms on the fiber. We do however need
to allow the hermitian form to take the value $\infty$ for some
vectors at some points in the base. We say that $|\xi|^2$ is a
quadratic form with values in $[0,\infty]$ on a vector space $V$ if
$|\xi|^2$ is a nonnegative quadratic form in the usual sense on a
subspace $V_0$ of $V$ and $|\xi|^2=\infty$ for any vector not in
$V_0$. Given a singular metric, $h$, the
norm function $|\xi|^2_h$ is a measurable function from
the total space of $E$ to $[0,\infty]$, whose restriction to any fiber
is a quadratic form in the sense just described.  

\smallskip

\noindent{\bf Definition 3.1} A singular metric  $h$ is {\it
  negatively curved}, if
$\log|\xi|^2_h$ is plurisubharmonic on the total space of $E$. 

Hence, a negatively curved metric is  upper
semicontinuous, in particular locally bounded from above, so it is a quadratic
form in the usual sense on any fiber. Note that the definition means
that a metric is 
negatively curved if and only if $\log|\xi(p)|^2$ is plurisubharmonic for any
choice of local holomorphic section of $E$. Since $\log(v_1+v_2)$ is
plurisubharmonic if $\log v_1$ and $\log v_2$ are plurisubharmonic, it
follows that a sum of 
negatively curved metrics is negatively curved.

Any choice of $h$ 
induces a dual metric $h^*$ on the dual bundle $E^*$. Then $h^*$ is
also an hermitian form on each fiber,  possibly taking the value
$\infty$ on some vectors. 

\noindent{\bf Definition 3.2} A singular metric $h$ is {\it positively
  curved}, if $h^*$ is negatively curved.

 Note that if  $h$ is smooth, then $h$ is positive if and only if 
$$
\Theta^h\geq 0
$$
in the sense of Griffiths.

\claim 3.1 Proposition|Let $E$ be a trivial holomorphic vector bundle
over a polydisk $U$, equipped with a singular metric denoted by $h$. If $h$ is
negatively curved over $U$, then  over any smaller polydisk
there is a sequence of smooth, strictly negatively curved metrics, $h_\nu$,
decreasing pointwise to $h$. Similarily, a
positively curved metric can be approximated from below by smooth strictly
positively curved metrics on smaller polydisks.
\endclaim

\proof. Assume $h$ is negatively curved. By adding a small multiple of a
trivial metric  ( this does not destroy the property of being
negatively curved) we may assume $h$ is strictly positive. Let
$$
h_\nu(p)=\int \chi_\nu(p-q)h(q),
$$
where $\chi_\nu$ is an approximate identity. 
Then
$$
|\xi|^2_{h_\nu}=\int \chi_\nu(p-q)|\xi|^2_h
$$
and it follows easily that $h_\nu$ is a smooth negatively curved metric. Then
multiply by $ 
e^{|p|^2/\nu}$ to  get strict positivity.
The corresponding statement for positively curved metrics follows by taking duals.
\hfill \qed

We will next give another criterion for positivity of a
bundle. Since the property is local we continue to assume that $E$ is
a trivial bundle over a polydisk, $U$. 
Given $(E,h)$ we  consider the bundle of unit balls
$$
\Omega=\{ (p,z)\in U\times\bC^N; |z|^2_{h(p)}\leq 1\}.
$$
The next preliminary statement will be very useful in what follows.

\claim 3.2 Proposition|Assume $h$ is smooth and strictly positively curved, and
let
$(p_0,z_0)$ be a point on $\partial\Omega$. Then there is a local
holomorphic section of $E$ near $p_0$, $s$, such that $s(p_0)=z_0$ and
$s$ takes values in $\Omega$.
\endclaim
\proof.
First choose $s$ holomorphic so that $s(p_0)=z_0$ and $D's=0$ at
$p_0$, where $D'$ 
is the $(1,0)$ part of the Chern connection of $h$. Then, if
$\Theta^h$ is the curvature of $h$,
$$
i\partial\bar\partial\log |s|^2=-\langle \Theta^h s,s\rangle_h/|s|^2_h<0.
$$
This means that there is a  second degree holomorphic polynomial,
$q(p)$, such that
$$
\log |s|^2_h -2\Re q
$$
has a local maximum at $p=p_0$. Then $s e^{-q}$ satisfies the
conclusion of the above proposition.
\hfill \qed
\medskip 
The main use of the positivity assumption for  bundles $E$ 
in this section and the next one will be to convert plurisubharmonic
functions on 
the total space of this bundle on plurisubharmonic functions on the base.
The precise statement is the following ( see also Demailly [12]):

\claim 3.3 Proposition|Let $V(p,z)$ be a continuous plurisubharmonic
function defined in some neighbourhood of $\Omega$ intersected with
the fiber $E_0$. Put
$$
V^*(p)=\sup_{z\in \Omega_p} V.
$$
If $h$ is positively curved, then $V^*$ is plurisubharmonic of $p$. Conversely,
if $V^*$ is plurisubharmonic of $p$ for any continuous
plurisubharmonic $V$, then $h$ is positively curved.
\endclaim
\proof.
We prove first the converse. Let $\xi(p)$ be a local holomorphic
section of $E^*$. Put
$$
V(p,z)=\log (|\xi(p)\cdot z|^2 +\epsilon).
$$
By assmption $V^*$ is plurisubharmonic. Letting $\epsilon$ go to 0, we
see that $\log|\xi|^2_{h^*}$ is plurisbharmonic. Hence $E^*$ is
negatively curved, so $E$ is positively curved.

Now assume that $E$ is positively curved, and let $h_\nu$ be sequence of smooth
strictly positively curved metrics increasing to $h$. Let $\Omega_\nu$ be the
corresponding bundles of unit balls. Then $\Omega_\nu$ decreases to
$\Omega$, so the corresponding functions
$$
V^*_\nu(p)=\sup_{z\in (\Omega_{\nu})_p} V
$$
decrease to $V^*$. Hence it is enough to prove our claim for the
metrics $h_\nu$; in other words we may assume that $h$ is smooth and strictly
positively curved. To prove that $V^*$ is plurisubharmonic we can take the base
dimension $m=1$, and moreover $V$ to be strictly plurisubharmonic. We
first claim that $V^*$ has no local maximum. For this, assume to get a
contradiction that 0 is a local maximum for $V^*$. Say
$$
V^*(0)=V(0,z_0).
$$
Take a local holomorphic section of $E$ near 0, satisfying the
conclusion of Proposition 2. Then $V(p,s(p))$ is strictly subharmonic of $p$ and
has a local maximum for $p=0$, which is impossible. 

Hence $V^*$ satisfies the maximum principle, and it follows that $V^*
+\Re q(p)$ also satisfies the maximum principle if $q$ is any
holomorphic polynomial. Thus $V^*$ is subharmonic.
\hfill \qed

\smallskip

Let us now return to the special case of vector bundles, $E$, that arise as
the direct image of $K_{X/Y}+L$, where $L$ is a line bundle over $X$,
fibered over $Y$.  Let
$$
B(x,x')=B_y(x,x')=\sum u_j(x)\otimes\bar u(x'),
$$
be the (relative) Bergman kernel defined at the beginning of section
1, but now considered also outside of the diagonal. The
plurisubharmonicity of $\log B(x,x)$ is equivalent to saying that if
$x(y)$ is a local holomorphic section of the fibration $X\mapsto Y$,
and if $a(y)$ is a local holomorphic section of $-(K_{X/Y}+L)$, then
$$
|a(y)|^2 B(x(y),x(y))
$$
is plurisubharmonic in $y$. Now choose for some $N>0$, $N$ holomorphic
sections of the fibration $x_j(y)$, and $N$ holomorphic sections of
$-(K_{X/Y}+L)$, $a_j(y)$, and consider the quadratic form
$$
\sum a_j\bar a_k B(x_j,x_k)=:B\langle a,a\rangle.
$$

\claim 3.4 Proposition| If $E$ is positively curved with respect to the
$L^2$-metric (1.1), then for any $N>0$ and any local sections as above
$$
B\langle a,a\rangle
$$
is a plurisubharmonic function of $y$.  Conversely, if this holds for
$N$ equal to the rank of 
$E$, then the $L^2$-metric of $E$ is positively curved.

\endclaim
\proof: Given local sections $x_j(y)$ and $a_j(y)$ as above, let
$$
\xi_y=\sum a_j(y) ev_{x_j(y)},
$$
with $ev_{x}$ denoting the evaluation functional at $x$. Then $\xi$ is
a holomorphic section of the dual bundle, and it is easily verified
that $B\langle a,a\rangle=|\xi|^2$. Hence $B\langle a,a\rangle$ is
plurisubharmonic if $E$ is positively curved. Conversely, if $N$ is the rank of
$E$ and if $x_j(y)$ are chosen
to be in general position for $y$ close to $y_0$, then any local
section of $E^*$ can be obtained in this manner. Hence $E^*$ is
negatively curved, so $E$ is positively curved.
\hfill \qed
\vskip 7pt

\noindent We state now the following generalization of the theorem 0.1, which will
play an important role in the next part of the paper. (It also motivates the
notions introduced and analyzed during this section.)

\claim Theorem 3.5| Let $p:X\mapsto Y$ be a smooth holomorphic
fibration, and let $L$ be a pseudoeffective line bundle over
$X$. Assume that
$$
E=p_*(K_{X/Y}+L)
$$
is locally free. Then $E$ is positively curved in the sense of Definition 3.2.
\endclaim
\proof. By the previous proposition we just need to check that
$B\langle a,a\rangle$ is plurisubharmonic for any choice of sections
$x_j$ and $a_j$. By the proof of that proposition, this means that
$|\xi|^2$ is plurisubharmonic. For the case of one single section
$x_j$ this follows from Theorem 0.1, as explained immediately before
the statement of proposition 3.4. The case of an $N$-tuple of sections
is proved in the same way: all the steps in the proof for $N=1$ work
in the same way for general $N$. 
\hfill \qed
\vskip 15pt
\section{\S 4 The Narasimhan-Simha metric}

\medskip
\noindent The Narasimhan-Simha metric is a metric on (possibly twisted)
multiples of the canonical bundle, $mK_X$, see e.g. Narasimhan-Simha [26],
Kawamata [18] and  Tsuji [32b] in a  more
general context.  
It is defined in a way
similar to the Bergman kernel metric, but using $L^p$- spaces instead
of $L^2$, with $p=2/m$. 

If $X$ is a compact complex manifold with a
holomorphic hermitian line bundle $(L,\phi)$, the Narasimhan-Simha 
metric (we will call it NS metric in the sequel) for
$mK_X+L$ is defined as the dual of the next metric on 
$-(mK_X+ L)$:
$$
\Vert \xi\Vert ^2: =\sup |\sigma(x)\cdot \xi|^2
$$
the $\sup$ being taken over all sections $\sigma$ to $mK_X +L$ such that
$$
\int_X |\sigma\wedge\bar\sigma|^{1/m}e^{-\phi/m}\leq 1
$$
and $\xi$ is some vector in 
$-(mK_X+ L)_x$.
\medskip

This defines $h_m$ as a metric on $mK_X +L$. The point of this
definition is that the quantity $|\sigma\wedge\bar\sigma|^{1/m}$
transforms as a measure on $X$ (in the case of $L$ trivial), 
so its integral can be defined without
choosing any metric on $K_X$. In the relative
situation, when the manifold $X$ is fibered over $Y$, we define
relative NS metrics in a similar way as in section 1, 
identifying $K_{X_y}$ with $K_{X/Y}$ restricted to $X_y$, and get a
metric on $mK_{X/Y}+L$. 

We shall now see that the methods of the previous section also give
positivity results for the relative NS metric. For this we will use a
family of $L^2$-metrics indexed by the sections to $mK_X +L$, defined
by
$$
\|\sigma\|^2_u=\int_X |\sigma|^2/|u|^{2-2/m} e^{-\phi/m}
$$
for sections $\sigma$ to $mK_X +L$. Thus the $L^p$-norm occuring in
the definition of the NS metric equals the 'diagonal' metric
$\|\sigma\|^2_\sigma $. This family of metrics  seems to be quite interesting
in their own right and play e.g. a fundamental role in Siu's proof of
the invariance of plurigenera. 

Let us now first consider an abstract model of the situation above: 
Take $E\rightarrow U$ to be a trivial vector bundle of rank $N$, with
an hermitian 
metric $h$, positive over 
$$
U=\{ (t,w); t\in \bC, |t|<1, w\in \bC^N\}.
$$
Choosing a trivialization, we can take
$$
E=U\times \bC^N=\{ (t,w,z); (t,w)\in U, z\in \bC^N\}.
$$
Let
$$
\Delta =\{ (t, w, z); w=z\},
$$
and
$$
\Omega_\Delta=\Omega\cap\Delta=\{(t,z,z); \|z\|^2_{h(z)}\leq 1\}.
$$
\vskip 7pt
We state now the next result, analog to the proposition 5.3 of the preceding section;
it has a major r\^ole in our investigation of the convexity properties of the NS metric.

\claim 4.1 Proposition|Assume $V(t,w,z)$ is continuous and
plurisubharmonic on $E$. Put
$$
V^\circ(t)=\sup_{\Omega_{\Delta, t}}V,
$$
where
$$
\Omega_{\Delta, t}=\{ (t,z,z); |z|^2_{h(z)}\leq 1\}.
$$
Then $V^\circ$ is subharmonic.
\endclaim
\proof.
Again, we may assume $h$ is smooth, strictly positively curved, $V$ strictly
plurisubharmonic, and it suffices to prove that $V^\circ$ has no local
maximum. Assume to the contrary that 0 is a local maximum point for
$V^\circ$, and say
$$
V^\circ(0)=V(0, z_0,z_0).
$$
Let $s$ be the local holomorphic section from Proposition 5.2, and put
$$
C=graf(s)\cap\Delta.
$$
Since the graf of $s$ has dimension $N+1$ and $\Delta$ has codimension
$N$, the dimension of $C$ is at least 1. Moreover
$$
(0,z_0,z_0)\in C\subset \Omega_\Delta.
$$
Thus $V$ has a local maximum on $C$, which contradicts the strict
plurisubharmonicity assumption.
\hfill \qed 
\vskip 10pt
\noindent As a consequence of the previous considerations, we obtain the next statement
which is a "twisted version" of the theorem 1 in the paper [18] by
Y. Kawamata. In this corollary we assume as a hypothesis that the
direct image of $mK_{X/Y}+L$ under $p$ is locally free, i e that for
any fiber of $p$,  any
holomorphic section of $(mK_{X/Y}+L)|{X_y}$ over the fiber, extends
holomorphically to some neighbourhood of the fiber. That this
hypothesis is satisfied for $L$ trivial (or $L$ having a semipositive
smooth metric) is exactly the invariance of (twisted) 
plurigenera.  
\smallskip
\claim 4.2 Corollary|Let $p:X\to Y$ be a smooth projective fibration and 
let $L$ be a holomorphic line bundle over $X$ endowed
with a  semipositively curved, possibly singular, metric. Assume that
the direct image of $mK_{X/Y}+L$ under the map $p$ is locally free. Then the
relative (and 
twisted) NS metric on $mK_{X/Y}+L$ has semipositive curvature current. 
\endclaim
\proof.
We will assume that in fact $X$ is fibered over the unit disk (see the
comments  
before the proof of the theorem 2.1 concerning the general case) and
let $F$ be the direct image 
of $mK_{X/Y}+L$. Since $F$ is locally trivial  it can be thought of as 
$$
F= D\times \bC^N.
$$
Now let $E$ be this bundle, pulled back over the total space of $F$
under the projection of $F$ to $D$. Then 
$$
E=D\times \bC^N_w\times\bC^N_z,
$$
and we consider $E$ as a bundle over $U:= D\times \bC^N_w$. Now, a
point $u=(t,w)$ in $U$ is a vector in $F$, ie a global section of $mK_{X_t}+L$
over $X_t$. Then $|u|^{2-2/m}e^{-\phi/m}$ is a metric on
$(m-1)K_{X_t}+L$ over $X_t$ and 
induces a hermitian metric on $E_u$. Since $u$ is holomorphic as a
function of $t$ and $w$ the metric  $|u|^{2-2/m}e^{-\phi/m}$ is pseudoeffective
on $(m-1)K_{X_t}+L$ considered as a bundle over $X\times
\bC^N_w$. Hence, by Theorem 3.5, the
induced hermitian metric on $E$ is positively curved . We can then apply
Proposition 4.1 to 
$$
V_\epsilon(t,w,z)=\log (|s(p)\cdot z|^2+\epsilon),
$$
with $s$ now equal to a section of $E^*$ giving point evaluations.
By Proposition 4.1 $V_\epsilon^\circ$ is plurisubharmonic. Hence
$V^\circ:=\lim V_\epsilon^\circ$ as $\epsilon$ goes to zero is also plurisubharmonic.
But  $V^\circ$ is precisely the NS metric. 
\hfill\qed

\vskip 7pt
\noindent {\bf Remark.} In the fundamental article [18] of Y. Kawamata, 
the convexity properties of the NS
metric are derived via the cyclic cover trick and subtle results in Hodge theory. 
In our proof we replace the cyclic covers by twisting the base 
with the space of all global holomorphic sections. More precisely, we have considered 
the product of the base $\Delta$ with the space of sections 
$\displaystyle H^0\big(X, mK_{X/\Delta}+ L\big)$; next, we have 
used the additional parameters to define
a metric on the bundle we are interested in. Moreover, one can see that given a map 
$p: X\to Y$ the above construction has no global meaning, but nevertheless we can use 
this local construction in order to derive  (global!) properties of
the NS metric,  
as indicated in the corollary above. Let us also mention that a more general version of the above statement will be given in the paper [4]: we will work in the context of the singular metrics and
surjective projective maps.

\vskip 10pt

Finally we mention one more consequence of the argument we just used,
where for simplicity we take $L$ to be tivial. 
Let $X$ be one fixed projective manifold and consider the line bundle
$m_1K_X$ over $X$. Let $m_0$ be a positive integer and let $u$ be a global holomorphic
section of $m_0K_X$; we use $u$ to define a metric on $(m_1-1)K_X$ as follows 
$$
|\sigma|^2_u =|\sigma|^2/|u|^{2(m_1-1)/m_0}
$$
and corresponding $L^2$-metrics
$$
\|\sigma\|^2_u :=\int_X |\sigma|^2_u 
$$
on $H^0(X,m_1K_X)$. Let $B_u(z)$ be the corresponding Bergman kernel.
\claim 4.3 Proposition|Let $p: \Gamma(X,m_0K_X)\times X\to X$
be the projection on the second factor and consider the bundle
$p^\star (m_1K_X)$. 
Then the Bergman kernel $B_u(z)$ defines a pseudoeffective metric on
$p^\star (m_1K_X)$. (Intuitively: it depends plurisubharmonically on
$u$ and $z$ together.)
In particular,  
let 
$$
Z=\{(u,z)\in \Gamma(X,m_0K_X)\times X; B_u(z)=0\}.
$$
Then either $Z$ is all of $ \Gamma(X,m_0K_X)\times X$ or $Z$ is
pluripolar in $ \Gamma(X,m_0K_X)\times X$
\endclaim
\proof.
Let $p$ be the projection map from $ \Gamma(X,m_0K_X)\times X$ to $
\Gamma(X,m_0K_X)$. This is a trivial fibration with fiber $X$ over 
$ \Gamma(X,m_0K_X)$ and its relative canonical bundle is equal to
the inverse image of $K_X$. The norms
$$
|\sigma|^2_u =|\sigma|^2/|u|^{2(m_1-1)/m_0}
$$
define a pseudoeffective metric on $(m_1-1)K_X$, so by Theorem 0.1 the
Bergman kernel (metric) is either identically 0 or pseudoeffective. In
the latter case it can vanish only on a pluripolar set, therefore the theorem is proved.

\vskip 30pt
\section{\S5 Proof of the corollary 0.2}   
\vskip 5pt

Before going into the details of the proof, we would like to mention the next 
result of  F. Campana (see [6]); it is a generalization of E. Viehweg--type semi-positivity theorems
(see e.g. [36]; in connection with this, see also the paper [19] of Y. Kawamata)
\claim Theorem ([6])|Let $p:X\to Y$ be a projective surjective map with connected fibers
between non-singular manifolds $X$ and $Y$. Let  
$L$ be a line bundle on $X$, whose first Chern class contains an effective
$\bQ$--divisor, which has trivial multiplier ideal sheaf when restricted to the general fiber of $p$. 
Then $p_*\bigl(m(K_{X/Y}+ L)\bigr)$ is weakly positive, for any positive and enough divisible integer $m$.
\endclaim

We do not intend to explain here the notion of {\sl weak positivity} introduced by 
E.~Viehweg in [37], [38]; however, let us point out a consequence of the previous result
which is stated in the paper [9].

\claim Theorem ([9])| Under the above assumptions, for any $\bQ$--ample divisor $A_Y\to Y$ some large enough multiple of the line bundle $K_{X/Y}+ L+ p^*A_Y$ is effective.
\endclaim

Therefore, our corollary 0.2 can be seen as a generalization and/or 
effective version of the result above (since there is no hypothesis on the multiplier sheaf in our statement and moreover we are able to control the restriction over some generic fiber of the sections produced by the theorem above).
\vskip 10pt
\noindent Next, we are going to explain the proof of 0.2. As we have already mentioned, our
arguments rely heavily on the theorem 0.1, and on the next version of the
Ohsawa-Takegoshi theorem, proved by L. Manivel (see [13], [23], [27] and the 
references quoted there).

\claim 5.1 Theorem ([23])|Let $X$ be a projective or Stein
$n$-dimensional 
manifold, and let $Z\subset X$ be the zero set of a holomorphic section
$s\in H^0(X, E)$ of a vector bundle $E\to X$; the subset $Z$ is assumed to be
non-singular, of codimension $r=\rank (E)$. 
Let $(F, h)$ be a line bundle on $X$,
endowed with a (possibly singular) metric $h$, such that:  

\smallskip
\item {(1)} $ \Theta_h(F)+ \sqrt {-1}\ddbar \log\vert s\vert ^2\geq 0$ as
current on $X$;
\smallskip
\item {(2)} $\displaystyle \Theta_h(F)+ \sqrt {-1}\ddbar \log\vert
s\vert ^2 \geq  
1/\alpha
{{\langle \Theta (E)s, s\rangle}\over {\vert s\vert ^2}}$ 
for some $\alpha\geq
1$;
\smallskip
\item {(3)} $\vert s\vert ^2\leq \exp (-\alpha)$ on $X$, and the restriction
of the metric $h$ on $Z$ is well defined.
\smallskip

\noindent Then every section $u\in H^0\bigl(Z, (K_X +
F_{\vert Z})\otimes {\cal I}(h_{|Z})\bigr)$ admits an extension
$\wt u$ to $X$ such that 

$$\int_X{{\wt u\wedge \ol {\wt u} \exp (-\varphi_F)}\over {\vert s\vert ^{2r}(-\log\vert
s\vert )^2}}\leq C_\alpha\int_Z{{|u|^2 \exp (-\varphi_F)dV_Z}\over 
{\vert \wedge ^r(ds)\vert ^2}},$$
where $h=e^{-\varphi_F}$, provided the right hand side is finite.
\endclaim 

\medskip

Notice that the integrand in the integral over $Z$ here is independent
of the metric chosen.

We will often use the Theorem 5.1 in the following situation: $(E,h_E)$ is a
holomorphic  line bundle with a smooth hermitian metric (remark that we do not require any curvature conditions for it) and $Z$ is a
hypersurface defined by a section $s$ to $E$. $\displaystyle (F',h_{F'})$ is 
another holomorphic hermitian line bundle over $X$ with a possibly
singular metric with semipositive curvature current, such that
$$
\Theta_h(F')\geq \epsilon\Theta_h(E)
$$
in the sense of currents. 
Then $F:= E+ F'$ endowed with the product metric $h_Eh_{F'}$ satisfies the curvature assumptions
in Theorem 5.1 if $\alpha$ is large enough. Hence any section of
$$
(K_X + F)|_Z
$$
that is square integrable over $Z$ extends to a section of $K_X+ F$
over all of $X$. Note also that by adjunction $(K_X+ E)|_Z=K_Z$, so we
are extending sections to $K_Z+ F'$, and
that the integrability assumption means that
$$
\int_Z c_nu\wedge\bar u \exp(-\varphi_F')<\infty.
$$

\noindent Also, in the papers quoted above,
the metric $h$ is non-singular; however, since the manifold 
$X$ above is projective, one can  derive
5.1 from the smooth case, by a regularization argument.
We refer to [34], [35] and [24] for the details concerning the
regularization process. 
Remark also that the only additional conditions needed in order to
allow singular metrics are the integrability, and the generic finiteness of
the restriction to the manifold we want to extend our form.

\vskip 10pt

The bundle we are interested in can be decomposed as follows:
$$m(K_{X/Y}+ L)
+ p^*A= K_X+ (m- 1)(K_{X/Y}+ L)+ L+ p^*(A- K_Y)$$
(the bundle $A$ will be chosen in a moment).
Our goal is to show that it is effective by extending (a multiple of)
the section $\displaystyle 0\neq u\in H^0\bigl(
X_y, (K_{X_y}+ L)\otimes {\cal I}(h)
\bigr)$ whose existence is given by the hypothesis. 

We denote by $h_B$ the metric on $K_{X/Y}+ L$ given by 
the theorem 1; the fact which will be crucial in what follows
is that for any positive integer $k$, we have 
$$\displaystyle 
\vert u^k\vert^2_{h_B^{\otimes k}}\leq {\cal O}(1)$$ 
uniformly on $X_y$, just by the construction
of the metric $h_B$.

\noindent We choose now the bundle $A_Y\to Y$
positive enough such that:   
\smallskip
\item {(a)}  $H^0(Y, A_Y)\neq 0$;
\smallskip
\item {(b)} The point $y\in Y$ is the common zero set of 
the sections $(s_j)$ of a ample line bundle $B\to Y$ and
$A_Y-K_Y\geq 2B$, in the sense that the difference is an ample line bundle.
\vskip 5pt
\noindent The bundle $(m- 1)(K_{X/Y}+ L)+ L+ p^*(A_Y- K_Y)$ is endowed with the
metric $$ h_m:= h_B^{\otimes (m-1)}\otimes h_L\otimes h_{A_Y-K_Y};$$
its curvature is semi-positive on $X$, and the restriction to
$X_y$ is well defined. The section we want to extend is
$\displaystyle u_m:= ~u^{\otimes m}\otimes s_{p^*A_Y}$, 
where $s_{p^*A_Y}$ is the pull-back of some 
non-zero section given by the property (a) above. 
By the property (b) the positivity conditions 
in the extension theorem of Ohsawa-Takegoshi are
satisfied with the bundle $F$ given by 
$$
F=(m- 1)(K_{X/Y}+ L)+ L+ p^*(A_Y- K_Y)
$$ 
since
$$\Theta _{h_m}\bigl( (m- 1)(K_{X/Y}+ L)+ L+ p^*(A_Y- K_Y)\bigr)+
\sqrt {-1}\ddbar \log\vert s\vert ^2\geq p^*\Theta (A_Y- K_Y-B).
$$
The right hand side  semi-positive and it
dominates  
the bundle $B$; thus the requirements (1), (2) and (3) are
verified. 

Now the integrability condition is obviously ok, since we have
$$\int_{X_y}\vert u_m\vert ^2_{h_m}dV\leq C\int_{X_y}\vert u\vert^2_{h_L}dV< \infty$$
(the first inequality comes from the fact that the section
$u$ is {\sl pointwise bounded} with respect to the fiberwise $L^2$ metric $h_B$,
and the second inequality is just the hypothesis). 
\vskip 10pt
Therefore, by the Ohsawa-Takegoshi theorem we can extend
the section $u_m$ over $X$ and the corollary 0.2 is proved.\hfill \qed

\claim Remark|{\rm In fact, we have just proved more than the corollary 0.2 states.
Indeed, we have shown that given $\displaystyle u \in H^0\bigl(
X_y, (K_{X_y}+ L)\otimes {\cal I}(h)
\bigr)$, then the section $u^{\otimes m}\otimes s_{p^*A_Y}$
extends over all of $X$ for any positive integer $m$.}
\endclaim
\vskip 25pt
\section{\S6 Maximal centers and asymptotic extensions}   
\vskip 10pt

\noindent In this paragraph we  give a 
complete proof of the theorem 0.3. We have already 
mentioned that this result is not new, but we feel that
the arguments 
which will be used here (which are mainly analytic)
could be very useful in other
problems in algebraic geometry.

To start with, we would 
like to recall the notions of {\sl restricted volume} and  
{\sl maximal center} of a line bundle, respectively of a $\bQ$-divisor.

\medskip
 
 \noindent {\bf Definition.} {\sl Let $E\to X$ be a line bundle, and let 
$V\subset X$ be an $d$-dimensional subset (which could be singular, but
reduced and irreducible). We define
$$H^0(X\vert V, mE):= \Im (H^0(X, mE)\to H^0(V, mE_V)),$$
and let $h^0(X\vert V,mE)$ be the dimension of this space.
Then the restricted volume of $E$ to $V$ is
$$\Vol_V(E):= \lim\sup {{d!}\over {m^d}} h^0(X\vert V, mL).$$}

\medskip 

\noindent If the line bundle $E$ is numerically effective, then
we have $\Vol_V(E)= E^d\cdot V$; this is a consequence of the Riemann-Roch theorem
(see [5] and also [14]). This notion turned out
to be very useful to deal with the linear systems which do have base-points, 
and one of the most spectacular result in the theory is the 
theorem 0.3 which will be discussed  next.
For further properties and consequences, we refer to [14].

\medskip

\noindent Let $D$ be an effective $\bQ$-divisor on $X$. We denote by 
${\cal I}(D)$ the multiplier ideal sheaf associated to $D$ (see e.g. [11], [25]).
To introduce the notion of maximal center of the pair 
$(X, D)$, let us assume that we have the decomposition
$D= S+ D^\prime$, where $S$ is a smooth hypersurface, and $D^\prime$
is an effective $\bQ$-divisor which does not contain $S$ in its support.
Then the adjunction formula gives $\displaystyle K_X+ D_{\vert S}
= K_S+ D^\prime_{\vert S}$, in particular the difference
$K_X+ D_{\vert S}- K_S$ is an effective $\bQ$-divisor on $S$.
As we will see in a moment, a maximal center is a substitute
in codimension greater than 1 of the hypersurface $S$ above.
We recall the following notions.
  
\vskip 15pt
 
 \noindent {\bf Definition.} 
 
 \noindent $1.$ {\sl We say that the pair $(X, D)$ is $klt$ at $x\in X$ if 
 $\displaystyle {\cal I}(D)_x= {\cal O}_{X, x}$; the support of the quotient
 $\displaystyle {\cal O}_{X}/{\cal I}(D)$ is denoted by $\Nklt(X, D)$.}
 
 \noindent $2.$ {\sl The log-canonical threshold of the pair $(X, D)$ at
 $x$ is the rational number 
 $$\lct_x(X, D):= \sup \{ t\in \bQ_+/{\cal I}(tD)_x= {\cal O}_{X, x} \} $$
 
 \noindent $3.$ We denote by ${\cal MC}(X, D)$ the set of 
 irreducible components $V$ of $\Nklt(X, D)$ such that the  log-canonical threshold of $(X, D)$ at the generic point 
of $V$ is equal to 1. The elements of this
 set are called maximal centers of $(X, D)$.}
 
 \medskip

\noindent We refer to the references [11], [19], [21] for a more extensive 
presentation and some properties of the notions introduced above;
here we will content ourself to recall some standard facts which will be
needed in what follows.

\smallskip
\noindent $\bullet $ Given a point $x\in X$, by a slight perturbation of the divisor $D$ we can
assume that there is at most one element of the set ${\cal MC}(X, D)$ containing $x$.

\smallskip
\noindent $\bullet $ If the maximal center $V$ is a divisor, then
we have the decomposition $D= V+ D^\prime$, where $D^\prime$ is an effective
$\bQ$-divisor, such that $V$ is not contained in the support of $D^\prime$.

\smallskip
\noindent $\bullet $ Assume now that the codimension of $V$ is 
at least 2, and assume that there exist an ample $\bQ$-line bundle
$A\to X$ such that $L:= A+ D$ is a genuine line bundle on $X$.
In this context we have the next ``concentration'' lemma of
Kawamata-Shokurov, see e.g. [30].

\claim 6.1 Lemma ([30])|Let $L= A+ D$ be a line bundle,
where $A$ and $D$ are $\bQ$-divisors ample and effective, respectively.
We assume that $V\in {\cal MC}(X, D)$ is of codimension at least 2.
Then there exist (another) decomposition $L= A+ \wt D$ such that: 

\noindent
\item {(A)} $V\in {\cal MC}(X, \wt D)$. 
\smallskip
\item{(B)}  There exist $\mu: Y\to X$ a sequence of blow-ups, such that we have
$$\mu^*(K_X+ L)= K_Y+ S+ \sum_j\nu_jE_j\leqno (6.1)$$
where $S$ and $E_j$ are in normal crossings, and the next properties hold: 
\smallskip
\item {(B.1)} \hskip 6pt $S$ is a smooth, $\mu$-contractible hypersurface of
$Y$, such that $\mu(S)= V$;  
\smallskip
\item {(B.2)} \hskip 6pt If for some index $j$ we have $V\subset \mu(E_j)$,
then the corresponding multiplicity $\nu_j$ is strictly less than 1.
\smallskip
\item {(B.3)} \hskip 6pt If the multiplicity $\nu_j$ is a negative number,
then the corresponding divisor $E_j$ is contracted by $\mu$.
\smallskip
\item {(B.4)} \hskip 6pt Some of the $E_j$ is ample on $Y$, and the
corresponding coefficient is positive and strictly less than 1.
\endclaim

\noindent Remark that if we don't need the property that a {\sl unique}
exceptional divisor with multiplicity 1 in the decomposition above 
dominates $V$, all the other facts in the lemma are immediate consequences of the 
fact that $V$ is the maximal center of $(X, D)$. However, this property will
simplify the arguments of the proof of the next corollary.
\vskip 10pt

\noindent We have divided the proof of 0.3 into two parts: the first
one uses techniques developed in connection 
with the invariance of plurigenera, and the second one uses the
semi-positivity of the relative canonical bundles 
(corollary 0.2).
 \vskip 15pt
 \subsection{\S6.1 Invariance of Plurigenera techniques}
 
 \medskip
The next result is due to S. Takayama in [30]; his proof follow very closely
the techniques introduced by Y.-T. Siu in [34].
Rather than reproducing here his arguments, we will use 
a more direct approach taken from [28] (see also the papers
[8], [30] and [36] for other variations and further results on this theme).
The version which is best adapted for our proof of 0.3 is the following.

\claim 6.1.1 Theorem ([30], [36])|Let $Z\subset X$ 
be a non-singular hypersurface 
and let $L= A+ D$ be a big line bundle, where $A$ and $D$ are 
$\bQ$-line bundles ample and effective, respectively. Let $A_X$ be a line bundle on $X$,
endowed with a smooth metric $\displaystyle h_{A_X}$ whose curvature is semi-positive.
Assume that $D$ does not contain $Z$ in its support
and that
${\cal I}(D_{\vert Z}) = {\cal O}_V$. Then for all $m\in \bZ_+$,
the restriction morphism  
$$H^0\bigl(X, m(K_X+ Z+ L)+ A_X\bigr)\to H^0\bigl(Z, (m(K_X+ Z+ L)+A_X)|_Z\bigr)$$
is surjective.
\endclaim

\medskip

\noindent The proof of this theorem is given in the last subsection in
order not to break the general argument.

\vskip 15pt

\subsection{\S6.2 End of the proof of 0.3}

\noindent Let $\mu: Y\to X$ be a log-resolution of the pair
$(X, D)$, such that the properties in the lemma 6.1 are satisfied.
We can re-write the relation (6.1) as follows

$$\mu^*(K_X+ L)= K_{Y}+ S+ D_Y- E_Y\leqno (6.2.1)$$
where $D_Y, E_Y$ are effective divisors
and $E_Y$ is contracted by $\mu$. 

If we denote by $D_Y^h$ the ``horizontal'' part of the divisor $D_Y$, 
(i.e. the components of $D_Y$
whose projection on $X$ contain $V$) 
then by the property $(B_2)$ of the lemma, we have
$[D_Y^h]= 0$ (we recall that the integral part, respectively the fractional part of a divisor 
is obtained by considering the integral part, respectively the fractional part
of its coefficients). 

We restrict the above equality to $S$ and use adjunction to get
$$\mu^*(K_X+ L)_{|S}= K_S+  D_Y- E_Y. $$
Next we subtract the inverse image of $K_V$ and split the $\bQ$--divisors 
$D_Y$ and $E_Y$ according to their integral/fractional parts; the above relation becomes
$$\mu^*(K_X+ L- K_V)_{\vert S}- [D_Y]- [-E_Y]_{\vert S}
= K_{S/V}+ \{D_Y\}+  \{-E_Y\}_{\vert S}.
\leqno (6.2.2)$$

\vskip7pt
\noindent We consider the $\bQ$-bundle $G:= \{D_Y\}+  \{-E_Y\}_{\vert S}$.
According with the notations in the lemma 6.1, we have 
$$G=  \sum_j(\nu_j- [\nu_j])E_j$$
and we can endow it with the canonical singular metric (observe in the first place that
$G$ is a genuine line bundle, by the relation (6.2.2) above). Along the next few lines, we
will highlight some of its metric properties, needed in what will follow.

\item {$(G_1)$}{\sl The multiplier ideal sheaf of the metric is trivial.}
Indeed, this is an obvious consequence of the ``normal crossings'' condition and the fact that the coefficients $G$
are in $(0, 1)$ (since the divisors $E_j$ are disjoint). 
\smallskip
\item {$(G_2)$}{\sl Given $v\in V$ generic, we have $\displaystyle H^0(S_v, K_{S_v}+ G)\neq 0$.}Indeed, by the
relation (6.~2.~2) and the property $(B_2)$ of the lemma 6.1, the bundle $\displaystyle K_{S_v}+ G_{\vert S_v}$ is isomorphic
with $\displaystyle -[-E_Y]_{S_v}$, so our claim is that this line bundle has at least a
non-zero section. It is the case since $\displaystyle -[-E_Y]$ has a section 
whose restriction to $S$ is non-trivial (recall the property (B.2)), and so will be the further restriction of this section of some generic fiber.
\smallskip
\item {$(G_3)$}{\sl The line bundle $G$ is big.} To see this, we recall that the property 
$(B_4)$ says that some of the $E_j$ in the decomposition (1) is ample on $S$ and that the corresponding coefficient $\nu_j$ belongs to $(0,1)$; therefore, the ample part is not lost when
we take $\nu_j- [\nu_j]$ in the construction of $G$.

\medskip

In conclusion, we have shown that we are in good position to apply our corollary 0.2: for any 
positive integer $m$, 
we get a non-trivial section  
 $\psi_m\in H^0\bigl(S, \mu^*A_X+ m(K_{S/V}+ G)_{\vert S})\bigr)$, where $A_X$ denotes here 
 some positive enough but fixed line bundle on $X$. 
The multiplication with 
 $\psi_m$ will define an 
injective map
$$\Psi_m: H^0(V, mK_V)\to H^0\bigl(S, \mu^*A_X+ m(K_Y+ S+ G)_{\vert
S}\bigr) \leqno (6.2.3)$$  
as follows: $\Psi_m(\tau)= \mu^*\tau\otimes \psi_m$.
\vskip 7pt
\noindent The result proved in the preceding paragraph 6.1 come now into the picture: the sections in the image of the morphism
$\Psi_m$ can be extended to $Y$ thanks to the theorem 6.1.1.
Indeed, the line bundle $G$ is big (see $(G_3)$ above) and the multiplier ideal of 
the restriction of its effective part to $S$ is trivial (by the same argument as in $(G_1)$).
\medskip

\noindent 
Let un consider now a section $u \in H^0(V, mK_V)$. By the discussion above, there exists a section
$\displaystyle \wt u\in H^0\bigl(Y, \mu^*A_X+ m(K_Y+ S+ G)\bigr)$ 
such that $\displaystyle \wt u_{\vert S}= \Psi_m\bigl(u \bigr)$. 
We denote by  $\displaystyle s_{m[D_Y]}$ the canonical section of the bundle
$m[D_Y]$, 
and then we have 
$\wt u\otimes s_{m[D_Y]}\in H^0\bigl(Y, \mu^*A_X+ m(\mu^*(K_X+ L)- [-E_Y])\bigr)$
by the definition of $G$ and the equality (6.2.1).

We recall now that every hypersurface in the support of the $\bQ$-divisor
$E_Y$ is contracted by $\mu$, therefore Hartogs type
arguments show that
$$\wt u\otimes s_{m[D_Y]}= \mu^*\bigl(J_m(u)\bigr)\otimes 
s_{-m[-E_Y]}\leqno (6.2.4)$$
for some section $J_m(u)\in H^0\bigl(X, A_X+ m(K_X+ L)\bigr)$.
In conclusion, we have defined an injective linear map
$$J_m:  H^0(V, mK_V)\to H^0\bigl(X\vert V, A_X+ m(K_X+ L)\bigr).$$

By hypothesis there exists a section $u\in H^0(X, m_0(K_X+ L)- A_X)$
such that $u_{\vert V}\neq 0$. Then by a multiplication with the section $u$, 
the above relation show the existence of a injective map
$$J_m^\prime:  H^0(V, mK_V)\to H^0\bigl(X\vert V, (m+ m_0)(K_X+ L)\bigr)$$
and this clearly implies the theorem 0.3, since the non-effective integer $m_0$ will disappear
as $m\to\infty$).

\claim Remark|{\rm As we have already mentioned in beginning the section 5, the version of the corollary 0.2 
needed above can also be obtained by the generalization due to F. Campana of some semi-positivity results of E. Viehweg. However, we  prefer to apply 
our result, not just because the technique needed to prove it is by far lighter, but also it provides 
an effective link between the canonical series of the manifolds $V$ and $X$, respectively
(see the comments page 15).}
\endclaim

\hfill \qed 

\medskip 

\subsection{\S6.3 Proof of Theorem 6.1.1}

\noindent We finally give the proof of Theorem 6.1.1.

 Let us fix an ample line bundle $B\to X$,
positive enough so that the next conditions hold true.

\noindent
\item {$\bullet$} For each $0 \leq p\leq m-1$, the bundle
$p(K_X+ Z+L)+ B$ is generated by its global sections, which we denote by
$(s^{(p)}_j)$;
\smallskip
\item {$\bullet\bullet$} Any section section of the bundle $m(K_X+Z+L)+ A_X+
B_{\vert Z}$ admits an extension to 
$X$.
\medskip
\noindent  Let $u\in H^0\bigl(Z, m(K_X+ Z+ L)+ {A_X}_{\vert Z}\bigr)$ be the section we want to extend. 
We consider the following statement:
\vskip 7pt
\noindent 
{$\bigl({\cal P}_{k,p}\bigr):$ 
\sl The sections $u^k\otimes s^{(p)}_j\in H^0\bigl(Z, (km+p)(K_X+ Z+ L)+ kA_X+ B_{\vert Z}\bigr)$
extend to $X$, for each $k\in \bZ_+$, $0 \leq p\leq m-1$ and $j= 1,..., N_p.$}
\vskip 7pt
\noindent If we can show that ${\cal P}_{k,p}$ is true for any $k$ and $p$, 
then another application of the extension theorem of Ohsawa--Takegoshi 
will end the proof; the argument goes as follows.

Let $\wt u^{(km)}_j\in H^0\bigl(X, km(K_X+ Z+ L)+ kA_X+ B\bigr)$
be an extension of $u^k\otimes s^{(0)}_j$.  We take $k\gg 0$, so that
$$ {{m-1}\over {mk}}B < A,
$$
 in the sense that the difference is
ample (we recall that $L= A+ D$). Let $h_A$ be a smooth, 
positively curved metric on $A-{{m-1}\over{mk}}B$. 
\vskip7pt
\noindent We apply now the extension theorem 5.1 with 
$F:= (m-1)(K_X+Z+ L) +L+ A_X+ Z;$
a metric $h$ on $F$ is constructed as follows:
\smallskip
\item {(1)} on the factor $(m-1)(K_X+Z+ L) +L+ A_X$ we take the
algebraic metric given by the family of sections
$\wt u^{(k)}_j$ (more precisely, we take the $\displaystyle {{m-1}\over
{mk}}$'th
root of this metric) multiplied with the canonical singular metric on $D$
twisted the metrics $h_A$ on $\displaystyle A- {{m-1}\over
{mk}}B$ and $\displaystyle h_{A_X}^{1/m}$ on $1/mA_X$ respectively; 
\smallskip
\item{(2)} we take an arbitrary, smooth metric $h_Z$
on the bundle associated to the divisor $Z$.
\smallskip
\noindent With our choice of the bundle $F$, 
the section $u$ we want to extend become a section of the adjoint 
bundle $K_X+ F_{\vert Z}$ and the positivity requirements (1) and (2) in the extension
theorem are satisfied since we have
$$ \Theta_h(F)+ \sqrt {-1}\ddbar \log\vert s\vert ^2_{h_Z}\geq \Theta_{h_A}(A- {{m-1}\over
{mk}}B)> 0$$
(in the notation above, $s$ is the canonical section of the bundle $\cO(Z)$).
Remark that the curvature of the bundle $\cO (Z)$ does not affect in any way
the $(1, 1)$--form above, since it cancel out by the Poincar\'e-Lelong identity. 
Also, by an appropriate choice of 
the constant $\alpha$ in the extension theorem and a rescaling of  the section
$s$, the conditions $(2)$ and $(3)$ will be satisfied as well.

Concerning the integrability of $u$,
remark that we have
$$\int_Z{{\vert u\vert^2\exp(-\varphi_D)}\over 
{\bigl(\sum_j\vert \wt u^{(km)}\vert^2\bigr)^{{m-1}\over {mk}}}}=
\int_Z{{\vert u\vert^2\exp(-\varphi_D)}\over 
{\bigl(\sum_j\vert u^{\otimes km}\otimes
s_j^{(0)}\vert^2\bigr)^{{m-1}\over {mk}}}}\leq
C\int_Z\exp(-\varphi_D)
$$
and this last integral converge by the hypothesis concerning the 
multiplier ideal sheaf of the restriction of $D$ to $Z$.

\noindent Thus, it is enough to check the property ${\cal P}_{k,p}$. 
For this, we will use an inductive procedure; the first steps are as follows.

\noindent (1)  For each $j= 1,...,N_0$, the section 
$u\otimes s^{(0)}_j\in H^0\bigl(Z, m(K_X+ Z+ L)+ A_X+ B_{\vert Z}\bigr)$
admits an extension $\wt u^{(m)}_j\in H^0\bigl(X, m(K_X+ Z+ L)+ A_X+ B\bigr)$,
by the property $\bullet\bullet$.

\noindent (2) We use the sections $(\wt u^{(m)}_j)$ to construct a metric
on $m(K_X+ Z+ L)+A_X+ B$; we multiply it with the (positively curved) metric of $L$
induced by $A$ et $D$; thus, we obtain a metric on 
the bundle $m(K_X+ Z+ L)+ L+ A_X+ B$.

\noindent (3) Let $F:= m(K_X+ Z+ L)+ L+ A_X+ B+ Z$; for each $j= 1,..., N_1$ let us consider the section 
$u\otimes s^{(1)}_j\in H^0\bigl(Z, (K_X+ F)_{\vert Z}\bigr)$.
It is integrable with respect to the metric produced at (2), since
$$\int_Z{{\Vert u\otimes s^{(1)}_j\Vert^2}\over {\sum_q \Vert u\otimes s^{(0)}_q\Vert^2}}
\exp (-\varphi_D)dV< \infty$$
again by hypothesis and $\bullet$.

\noindent (4) We apply the extension theorem (as above, we take an arbitrary smooth metric
on $\cO(Z)$) and we get
$(\wt u^{(m+ 1)}_j)\in H^0\bigl(X, (m+1)(K_X+ Z+ L)+A_X+ B\bigr)$ whose
restriction on $Z$ is precisely
$u\otimes s^{(1)}_j$.

\vskip 5pt Now the assertion  ${\cal P}_{k,p}$ will be obtained by iterating
the
procedure (1)-(4) several times. Indeed, assume that the 
proposition ${\cal P}_{k,p}$ has been checked, and consider the set of global sections 
$$(\wt u^{(km+p)}_j)\in H^0\bigl(X, (km+p)(K_X+ Z+ L)+ kA_X+ B\bigr)$$
which extend $u^k\otimes s^{(p)}_j$. The metric associated to them twisted with the metric of $L$ induce a metric on the bundle
$(km+p)(K_X+ Z+ L)+ L+ kA_X+ B$. Now if $p< m-1$, we can consider the family of 
sections $u^k\otimes s^{(p+1)}_j\in H^0(Z, K_X+ F_{\vert Z})$, where 
$F:=(km+p)(K_X+ Z+ L)+ L+ kA_X+ B+Z$. Each of them is integrable with respect to the
previous metric (by the same arguments as in (3) above) and the extension theorem 3.1 show that   
${\cal P}_{k,p+1}$ is verified. In the remaining case $p= m-1$, we will consider the 
bundle $F$ as above twisted with $A_X$ and its corresponding metric. Over $Z$, the section we will consider are 
$u^{k+1}\otimes s^{(0)}_j\in H^0(Z, K_X+ F)$; the rest of the proof
goes along the same lines, 
thus we skip it.
\hfill \qed
\vskip 5pt
\vfill \eject 
\section{References}

\bigskip

{\eightpoint

\bibitem [1]&Berndtsson, B.:&Positivity of direct image bundles and
  convexity on the space of Kahler metrics &math arXiv: math.CV/0608385&  

\bibitem [2]&Berndtsson, B.:&  Subharmonicity properties of the Bergman kernel and some other functions associated to pseudoconvex domains;& (math arXiv:mathCV 0505469, to appear in Ann Inst Fourier)&

\bibitem [3]&Berndtsson, B.:& Curvature of Vector bundles associated to holomorphic fibrations;& to appear in Ann.\ of Maths.\ (2007)&

\bibitem [4]&Berndtsson, B., P\u aun, M.:& Bergman kernels and the pseudo-effectivity of the relative canonical bundles --Part two--;& in preparation&

\bibitem [5]&Boucksom, S.:& On the volume of a big line bundle;&
  Intern.\ J.\ Math.\ {\bf 13} (2002), 1043--1063&
  
 \bibitem [6]&Campana, F.:& Special varieties, orbifolds and classification theory;&
  Ann.\ Inst.\ Fourier {\bf 54} (2004), 499--665&
  
 \bibitem [7]&de Cataldo, M.A.:& Singular hermitian metrics on vector bundles
           &J. Reine Angew. Math. {\bf 502} (1998), 93-122&

\bibitem [8]&Claudon, B.:& Invariance for multiples of the twisted canonical bundle& math.AG/0511736, to appear in Annales de l'Institut Fourier&

\bibitem [9]&Debarre, O.:& Syst\`emes pluricanoniques sur les vari\'et\'es de type g\'en\'eral;&
s\'eminaire Bourbaki {\bf 970}, 2006-2007&

\bibitem [10]&Demailly, J.-P.:& Singular hermitian metrics on positive line bundles;& Proceedings of the Bay\-reuth conference ÒComplex algebraic varietiesÓ, April 2- 
{\bf 6}, 1990, edited by K. Hulek, T. Peternell, M. Schneider, F. Schreyer, Lecture 
Notes in Math.\ {\bf 1507}, Springer-Verlag, 1992&

\bibitem [11]&Demailly, J.-P.:& A numerical criterion for very ample line bundles;& J.\  Differential Geom.\ {\bf 37} (1993), 323-374&

\bibitem [12]&Demailly, J-P.:& Pseudoconvex-concave duality and
  regularization of currents;&  Several complex variables (Berkeley,
  CA, 1995--1996),  233--271, Math. Sci. Res. Inst. Publ., 37,
  Cambridge Univ. Press, Cambridge, 1999.&

 \bibitem [13]&Demailly, J.-P.:&  On the Ohsawa-Takegoshi-Manivel  
extension theorem;& Proceedings of the Conference in honour of the 85th birthday of Pierre Lelong, 
Paris, September 1997, Progress in Mathematics, Birkauser, 1999&

 \bibitem [14]&Ein, L. {\it et al}:& Restricted volumes and base loci of linear series;&available at  math.AG/0607221&
  
 \bibitem [15]&Fujita, T.:& On K\"ahler fibre spaces over curves;& J. \ Math. \ Soc. \ Japan {\bf 30} (1978),  779Ð794&

 \bibitem [16]&C. D. Hacon, J. McKernan:& Boundedness of pluricanonical maps of varieties of general type;&
 Invent.\ Math.\ Volume {\bf 166}, Number 1 / October, 2006, 1-25&

 \bibitem [17]&Kawamata, Y.:& Minimal models and the Kodaira dimension of algebraic fibre
spaces;& Journ.\ Reine Angew.\ Math.\  {\bf 363} (1985) 1-46&
 
 \bibitem [18]&Kawamata, Y.:& Kodaira dimension of algebraic fiber spaces over curves;&  Invent.\ Math.\ {\bf 66} (1982)  no. 1, 57--71&

\bibitem [19]&Kawamata, Y.:& Subadjunction of log canonical divisors;& Amer.\ J.\ Math.\  {\bf 120} (1998) 893--899&

 \bibitem [20]&Koll\'ar, J.:& Subadditivity of the Kodaira dimension: Fibres of general type;& Algebraic geometry, Sendai 1985, Advanced Studies in Pure Math.\ {\bf 10} (1987)
361Ð398&
 
 \bibitem [21]&Lazarsfeld, R.:& Positivity in Algebraic Geometry;& Springer, Ergebnisse der Mathematik und ihrer Grenzgebiete&
 
 \bibitem [22]&Maitani F., Yamaguchi H.:& Variation of Bergman metrics on Riemann surfaces;&Mathematische Annalen
Volume {\bf 330}, Number 3 / November, 2004, 477-489&

 \bibitem [23]&Manivel, L.:& Un th\'eor\`eme de prolongement L2 de sections holomorphes d'un 
 fibr\'e hermitien;& Math.\ Zeitschrift {\bf 212} (1993), 107-122&

 \bibitem [24]&McNeal J., Varolin D.:&  Analytic inversion of adjunction: $L\sp
2$ extension theorems with gain;& Ann. Inst. Fourier (Grenoble) 57
(2007), no. 3, 703--718& 

\bibitem [25]&Nadel, A.M.:&Multiplier ideal sheaves and K\"ahler-Einstein 
  metrics of positive scalar curvature;& Ann.\ of Math.\ {\bf 132} (1990), 549Ð596&

\bibitem [26]&Narasimhan, M. S.; Simha, R. R.;&Manifolds with ample
  canonical class. ;&  Invent. Math.  {\bf 5}  (1968) 120--128. &

\bibitem [27]&Ohsawa, T., Takegoshi, K.\ :&  On the extension of $L^2$
holomorphic functions;& Math.\ Z.,
{\bf 195} (1987), 197--204&

\bibitem [28]&P\u aun, M.:& Siu's Invariance of Plurigenera: a One-Tower Proof ;&preprint 2005, to appear in J.\ Diff.\ Geom.\ &

\bibitem [29]&Takayama, S:& On the Invariance and Lower Semi--Continuity
of Plurigenera of Algebraic Varieties;& J. Algebraic Geom.  {\bf 16 } (2007), no. 1, 1--18&

\bibitem [30]&Takayama, S:& Pluricanonical systems on algebraic varieties of general type;& Invent.\ Math.\
Volume {\bf 165}, Number 3 / September, 2005, 551-587&

\bibitem [31]&Tsuji, H.:& Variation of Bergman kernels of adjoint line
  bundles;& math.CV/0511342 &
 
\bibitem [31b]&Tsuji, H.:& Dynamical construction of K\"ahler-Einstein
  metrics;& math.AG/0606023 &

\bibitem [32]&Tsuji, H.:&  Pluricanonical systems of projective
  varieties of general type I;& Osaka J of Math {\bf 43-4} (2006), 967-995&

\bibitem [32a]&Tsuji, H.:&  Pluricanonical systems of projective
  varieties of general type II;& Osaka J of Math  {\bf 44-3} (2007), 723-734&

\bibitem [32b]&Tsuji, H.:& Curvature semipositivity of relative
  pluricanonical systems;& math.AG/0703729&

\bibitem [33]&Siu, Y.-T.:& Invariance of Plurigenera;& Inv.\ Math.,
{\bf 134} (1998), 661-673&

\bibitem [34]&Siu, Y.-T.:& Extension of twisted pluricanonical
  sections with plurisubharmonic weight and invariance of
  semipositively twisted plurigenera for manifolds not necessarily of
  general type;& Complex geometry (G\"ottingen, {\bf 2000}),
  223--277, Springer, Berlin, 2002& 

\bibitem [35]&Siu, Y.-T.:& Multiplier ideal sheaves in complex and
  algebraic geometry;& Sci.\ China Ser.  {\bf A 48}, 2005& 

\bibitem [36]&Varolin, D.:&  A Takayama-type extension theorem;&
  math.CV/0607323, to appear in Comp.\ Math& 
     
\bibitem [37]&Viehweg, E.:& Weak positivity and the additivity of the
  Kodaira dimension for 
certain fibre spaces;& Algebraic Varieties and Analytic Varieties, Advanced
Studies in Pure Math.\ {\bf 1} (1983) 329Ð353&
     
\bibitem [38]&Viehweg, E.:& Quasi-Projective
Moduli for
Polarized Manifolds;& Springer-Verlag, Berlin, Heidelberg, New York, 1995
as: Ergebnisse der Mathematik und ihrer Grenzgebiete, 3. Folge, Band 30&
}

\bigskip
\noindent
(version of January 21, 2008, printed on \today)
\bigskip\bigskip
{\parindent=0cm
Bo Berndtsson,  
bob@math.chalmers.se\\
Mihai P\u aun,
paun@iecn.u-nancy.fr
}
\end